\documentclass[10pt]{amsart}

\usepackage[english,francais]{babel}
\usepackage[T1]{fontenc}
\usepackage{amsmath,amssymb}
\usepackage[all]{xy}
\usepackage{graphicx}

\newcommand{\lgw}{\longrightarrow}
\newcommand{\lgm}{\longmapsto}
\newcommand{\ovl}{\overline}

\newcommand{\Spec}{\text{Spec}\,}
\newcommand{\ord}{\text{ord}}
\newcommand{\Ker}{\text{Ker}}

\newcommand{\D}{\Delta}

\newcommand{\wdt}{\widetilde}

\newcommand{\Z}{\mathbb{Z}}

\renewcommand{\k}{\Bbbk}
\newcommand{\E}{\mathbb{E}}

\newcommand{\U}{\mathbb{U}}
\newcommand{\F}{\mathbb{F}}
\newcommand{\R}{\mathbb{R}}
\newcommand{\Chi}{\mathbf{\chi}}

\renewcommand{\H}{\mathbb{H}}
\newcommand{\K}{\mathbb{K}}
\newcommand{\DL}{\mathbb{L}}

\newcommand{\N}{\mathbb{N}}
\newcommand{\A}{\mathbb{A}}
\newcommand{\C}{\mathbb{C}}
\newcommand{\Q}{\mathbb{Q}}
\renewcommand{\a}{\alpha}
\renewcommand{\b}{\beta}
\newcommand{\g}{\gamma}

\newcommand{\e}{\varepsilon}

\newtheorem{propdef}{Proposition-D\'efinition}[section]
\newtheorem{lemm}[propdef]
{Lemme }
\newtheorem{theo}[propdef]
{Th\'eor\`eme }

\newtheorem{coro}[propdef]
{Corollaire }
\newtheorem{prop}[propdef]{Proposition}

\theoremstyle{definition}
\newtheorem{exem}[propdef]{Exemple }
\newtheorem{rema}{Remarque}
\newtheorem{defi}[propdef]{D\'efinition }

\author{Guillaume Rond}
\address{\newline
  Department of Mathematics 
\newline
University of Toronto
 \newline
Toronto, Ontario
  \newline
 Canada M5S 3G3}
\email{rond@picard.ups-tlse.fr}

\title[S\'eries de Poincar\'e motiviques]{S\'eries de Poincar\'e motiviques d'un germe d'hypersurface irr\'eductible quasi-ordinaire }

\begin{document}

\begin{abstract}
Nous donnons ici une description combinatoire, faisant intervenir les exposants caract\'eristiques de la singularit\'e, des arcs tronqu\'es trac\'es sur un germe irr\'eductible d'hypersurface quasi-ordinaire. Cela nous permet d'obtenir une expression inductive des s\'eries de Poincar\'e de ce type de singularit\'e.
\end{abstract}

\maketitle
\tableofcontents

\section{S\'eries de Poincar\'e motiviques}
Soit $(X,0)$ un germe de vari\'et\'e analytique sur un corps $\k$ de caract\'eristique nulle. Soit $p$ un entier naturel. Nous d\'efinissons l'espace des jets d'ordre $p$, not\'e $X_p$, comme \'etant la vari\'et\'e alg\'ebrique sur $\k$ dont les points $\K$-rationnels, pour toute extension de corps $\K$ de $\k$, sont les  $\K[t]/t^{p+1}$-points de $(X,0)$. C'est-\`a-dire que nous avons $X_p=\{\varphi : \Spec \k[[t]]/t^{p+1}\lgw (X,0)\}$. Dans le cas particulier o\`u $(X,0)$ est un germe de vari\'et\'e analytique d\'efinie par les \'equations $f_i(x)=0$, pour $i=1,...,\,r$ et $x=(x_1,...,\, x_s)$, alors  $X_p$ est la vari\'et\'e affine d\'efinie par les \'equations en les variables $x_{j,k}$ pour $k=1,...,\,p$ et $j=1,...,\,s$, provenant du fait que $f_i(x_{j,1}t+\cdots+x_{j,p}t^p)=0$ mod $t^{p+1}$ pour tout $i$.\\
La limite projective de ces vari\'et\'es, appel\'ee espace des arcs sur $X$, est not\'ee $X_{\infty}$ et n'est en g\'en\'eral pas de type fini sur $\k$.\\
Nous avons les morphismes naturels de troncations
$$\pi_p\ :\ X_{\infty}\lgw X_p \ \text{ et } \ \pi_{p,q}\ :\ X_p\lgw X_q\text{ pour $p\geq q$} \ .$$
Nous nous int\'eressons ici au comportement des arcs tronqu\'es, c'est-\`a-dire aux $\pi_p(X_{\infty})$ quand $p$ varie. Nous savons, d'apr\`es le th\'eor\`eme de Greenberg \cite{Gr}, que ce sont des ensembles constructibles. On peut donc consid\'erer leur image dans l'anneau de Grothendieck $K_0(Var_{\k})$ des vari\'et\'es sur $\k$ \cite{D-L1}. Plus pr\'ecis\'ement nous nous int\'eressons \`a la s\'erie de Poincar\'e g\'eom\'etrique $P_{geom,X,0}(T):=\sum_{p\geq 0}[\pi_p(X_{\infty})]T^p$ o\`u $[Y]$ repr\'esente la classe de la vari\'et\'e $Y$ dans $K_0(Var_{\k})$. Denef et Loeser ont montr\'e que cette s\'erie \'etait rationnelle avec un d\'enominateur qui s'\'ecrit sous forme du produit de termes de la forme $1-\DL^aT^b$ o\`u $\DL:=[\A^1_{\k}]$ et $a\in\Z$ et $b\in\N\backslash\{0\}$ (cf. \cite{D-L1}). Cependant la preuve utilise \`a la fois la r\'esolution des singularit\'es de la vari\'et\'e et un th\'eor\`eme d'\'elimination des quantificateurs d\^u \`a Pas \cite{Pa}, et n'apporte aucune information quantitative, en particulier sur les p\^oles. Cette s\'erie a jusqu'\`a pr\'esent \'et\'e calcul\'ee pour les branches planes (cf. \cite{D-L2}) et les singularit\'es de surfaces toriques normales (cf. \cite{L-J-R} et \cite{Ni1}).\\
Par ailleurs, on peut aussi consid\'erer $\varphi_p$,  la formule dans le langage de premier ordre de $\k[[T]]$, d\^u \`a Pas \cite{Pa}, qui d\'efinit $\pi_p(X_{\infty})$, qui est un ensemble constructible, et regarder sa mesure arithm\'etique $\Chi_c(\varphi_p)$ dans l'anneau de Grothendieck $K^v_0(Mot_{\k,\,\ovl{\Q}})_{\Q}$, c'est-\`a-dire l'anneau de Grothendieck des motifs de Chow sur $\k$ \`a coefficients dans $\ovl{\Q}$ tensoris\'e avec $\Q$ (cf. \cite{D-L2} ou \cite{Ha} pour une introduction). 
Une autre s\'erie int\'eressante est alors la s\'erie d\'efinie par $\sum_{p\geq 0}\Chi_c(\varphi_p)T^p$. Cette s\'erie se sp\'ecialise pour tout $q$ premier, sauf un nombre fini, en la s\'erie $\sum_{p\geq 0}N_{q^p}(X,0)T^p$, o\`u $N_{q^p}(X,0)$ est le cardinal des $\Z/q^p\Z$-points de $(X,0)$ qui se rel\`event en des $\Z_q$-points de $(X,\,0)$. Denef et Loeser ont montr\'e le m\^eme r\'esultat de rationnalit\'e pour cette s\'erie que pour la s\'erie g\'eom\'etrique (cf. \cite{D-L2}). Cette s\'erie a \'et\'e calcul\'ee pour les branches planes (cf. \cite{D-L2}) et les singularit\'es de surfaces toriques normales (cf. \cite{Ni1}). Pour les branches planes ces deux s\'eries diff\`erent et pour les surfaces normales toriques, J. Nicaise montre l'\'egalit\'e.\\
\\
Nous avons effectu\'e ici le calcul des s\'eries g\'eom\'etrique et arithm\'etique d'un germe d'hypersurface irr\'eductible quasi-ordinaire. Ce type de singularit\'e g\'en\'eralise les singularit\'es de courbes planes dans le sens o\`u il existe une param\'etrisation de ces singularit\'es \`a l'aide de s\'eries fractionnaires \`a plusieurs variables dont le d\'enominateur est born\'e. L'ensemble des exposants, apparaissant dans l'\'ecriture des ces s\'eries, appartient au groupe engendr\'e par un nombre fini d'exposants, appel\'es \textit{exposants caract\'eristiques}, qui g\'en\'eralisent les exposants caract\'eristiques d'une courbe plane.\\
Pour calculer la mesure  motivique de l'ensemble des arcs tronqu\'es \`a l'ordre $p$, nous d\'ecomposons cet ensemble en deux ensembles constructibles :  l'ensemble des arcs tronqu\'es qui ne se rel\`event pas en arcs inclus dans le compl\'ementaire du tore et son compl\'ementaire. Nous donnons d'abord une caract\'erisation combinatoire  des arcs tronqu\'es qui ne se rel\`event pas en arcs inclus dans le compl\'ementaire du tore.  Cela nous permet de calculer la mesure motivique de l'ensemble de ces arcs tronqu\'es. Enfin nous donnons une formule de r\'ecurrence sur la dimension du germe pour la mesure de son compl\'ementaire. Nous obtenons alors des formules g\'en\'erales, inductives sur la dimension de l'hypersurface, de ces deux s\'eries. Nous donnons enfin une formule explicite de ces s\'eries dans le cas o\`u les coordonn\'ees des exposants caract\'eristiques de la singularit\'e sont sup\'erieures \`a 1, c'est-\`a-dire quand la projection de celle-ci est ``tr\`es'' transverse (th\'eor\`eme \ref{exp}).\\
\\
Je tiens  \`a remercier ici  M. Lejeune-Jalabert pour avoir fait preuve de patience \`a l'\'ecoute de ces r\'esultats et pour ses pr\'ecieux commentaires. Je remercie aussi J. Nicaise pour m'avoir, le premier, parl\'e des s\'eries de Poincar\'e motivique lors du GAEL XII, et Helena Cobo Pablos et Pedro Gonzalez-Perez pour m'avoir indiqu\'e une erreur dans la premi\`ere version de ce travail.

\section{Singularit\'es quasi-ordinaires}
Nous rappelons ici la d\'efinition de singularit\'e  quasi-ordinaire et les propri\'et\'es de ces singularit\'es dont nous aurons besoin. 
\begin{defi}(cf. \cite{G-P1} par exemple) Soit $f\in\C\{X_1,...,\,X_m\}[Y]$  un polyn\^ome distingu\'e. On dit que $f$ est quasi-ordinaire si  son discriminant $\D_Y(f)$ a un terme dominant, c'est-\`a-dire si il s'\'ecrit $X^{\a}u$ o\`u $\a\in\Z_{\geq0}^{m}$ et $u$ est inversible. G\'eom\'etriquement, cela revient \`a dire que le discriminant de la projection du germe $(X,0)\subset\C^{m+1}\lgw \C^{m}$ qui envoie le point de coordonn\'ees $(x_1,...,\,x_m,\,y)$ sur le point de coordonn\'ees  $(x_1,...,\,x_m)$ est \`a croisements normaux.
\end{defi}
Nous avons alors la proposition

\begin{prop}\cite{Ab}
Les racines de $f$ vu comme polyn\^ome en $Y$ sont dans 
$$\wdt{\C\{X\}}=\lim_{\begin{array}{c}
\rightarrow\\
k\in\N^*\\
\end{array}}\C\{X^{\frac{1}{k}}\}.$$
De plus les racines de $f$ sont deux \`a deux comparables, c'est-\`a-dire que pour toutes $\xi$ et $\xi'$ racines de $f$, $\xi-\xi'$ s'\'ecrit $X^{\a}u$ o\`u $\a\in\Q_{\geq0}^{m}$ et $u$ est inversible.

\end{prop}

Le semi-groupe $\N^{m+1}$ \'etant muni de l'ordre lexicographique $\leq_{lex}$, nous pouvons d\'efinir l'ordre suivant sur $\N^m$ :\\
Si $(\a_1,...,\,\a_m)$ et $(\b_1,...,\,\b_m)\in\N^m$, $(\a_1,...,\,\a_m)\leq (\b_1,...,\,\b_m)$ si et seulement si $(|\a|,\a_1,...,\,\a_m)\leq_{lex}(|\b|,\b_1,...,\,\b_m)$.
 Ainsi $\C\{X_1,...,\,X_m\}$ est muni d'une valuation $\nu_X$, \`a valeurs dans $\N^{m}$ muni de l'ordre d\'efini pr\'ecedemment, qui \`a un mon\^ome $X^{\a}$ associe $(\a_1,...,\,\a_m)$.  La valuation $\nu_X$ s'\'etend alors \`a $\wdt{\C\{X\}}$ (mais son groupe de valeurs est alors $\Q^m$).\\

\begin{theo}\label{quasi}\cite{Ga},\cite{Li2}
Soit $f$  irr\'eductible et quasi-ordinaire. Alors nous avons :
\begin{enumerate}
\item Si $deg_Y(f)=n$ alors  $f$ a $n$ racines distinctes dans  $\C\{X^{\frac{1}{n}}\}$.
\item Si $\xi$ est une racine de $f$ dans $\C\{X^{\frac{1}{n}}\}$, alors il existe des \'elements de $(\frac{1}{n}\Z)^m$, strictement ordonn\'es, $a(1)<\cdots<a(g)$ (i.e. $a_k(1)\leq\cdots\leq a_k(g)$ pour tout $k$ et $a(i)\neq a(j)$ pour $i\neq j$)
 tels que l'on puisse \'ecrire
$$\xi=\xi_0+\xi_1+\cdots+\xi_g$$ 
$$\text{avec } \ \xi_0\in\C\{X\},$$
$$X^c \text{ appara\^it dans } \xi\Longrightarrow c\in \Z^m+\sum_{a(i)\leq c}a(i)\Z,$$
$$X^c \text{ appara\^it dans } \xi_i\Longrightarrow c\in \Z^m+\sum_{j\leq i}a(j)\Z,$$
$$\text{et }\ \nu_X(\xi_k)=a(k)\ \text{ pour tout }k.$$
\item Si $\xi$ est racine de $f$ dans $\C\{X^{\frac{1}{n}}\}$, alors l'ensemble des  racines de $f$ est l'ensemble form\'e des $\xi(w_1X_1^{\frac{1}{n}},...,\, w_mX_m^{\frac{1}{n}})$ o\`u les $w_k$ parcourent l'ensemble des racines $n$-i\`emes de l'unit\'e.\\
\end{enumerate}
\end{theo}

\begin{rema}
Quitte \`a faire le changement de variables 
$$X_k\lgm X_k,\ \forall k$$
$$\text{ et }\, Y\lgm Y+\xi_0$$
nous pouvons supposer que $\xi_0=0$ dans le th\'eor\`eme pr\'ec\'edent.\\

\end{rema}

\begin{defi}
Les $a(k)$ du lemme pr\'ec\'edent sont appel\'es  les exposants caract\'eristiques de $f$.
\end{defi}
Nous pouvons  aussi d\'efinir l'ensemble des exposants caract\'eristiques comme \'etant l'ensemble
$$\left\{\nu_X(\xi-\xi'),\,f(\xi)=f(\xi')=0,\,\xi\neq\xi'\right\}.$$\\
\\
Nous pouvons d\'efinir les r\'eseaux  $M=M_0:=\Z^m$ et $M_k:=\Z^m+\sum_{l\leq k}a(l)\Z^m$ pour $1\leq k\leq g$ et les r\'eseaux duaux $N_i:=\check{M_i}$ et $N=N_0$. Nous d\'efinissons les entiers caract\'eristiques $n_k$ du germe d'hypersurface comme \'etant les indices des $M_{k-1}$ dans $M_k$ :\\
$$n_k=[M_k:M_{k-1}].$$
Nous posons aussi $n_0=1$ et $n_{-1}=0$.
Nous notons
$$e_{k-1}=n_k...n_g\ \text{pour } k=1,...,g.$$
Notons $L=L_0$ le corps des fractions de $\C\{X\}$. Nous pouvons remarquer que pour $k=1,...,g$,
$$e_k=[L[\xi]:L[X^{a(1)},...,\,X^{a(k)}]]=[L[\xi]:L[\xi_1+\cdots+\xi_k]]$$
$$\text{ et }\, n_k=[L[X^{a(1)},...,\,X^{a_k}]:L[X^{a(1)},...,\,X^{a(k-1)}]]$$
$$=[L[\xi_1+\cdots+\xi_k]:L[\xi_1+\cdots+\xi_{k-1}]].$$
En particulier nous avons $e_0=n=n_1...n_g$.\\
Nous pouvons aussi d\'efinir les vecteurs $\g(k)$ par :
$$\g(1):=a(1)$$
$$\g(k+1):=n_k\g(k)+a(k+1)-a(k)$$
Dans le cas $m=1$, ce sont les $n\g_i$ sont des g\'en\'erateurs du semi-groupe de l'ensemble des multiplicit\'es d'intersection $(C,\,X)_0$ o\`u $C$ parcourt l'ensemble des germes en 0 de courbes planes non contenues dans $X$ (cf. \cite{Z}).\\
Nous avons le 
\begin{lemm}\cite{G-P2}
Le sous-r\'eseau de $M_g$ engendr\'e par $\g(1),...,\,\g(k)$ est \'egal \`a $M_k$ et l'ordre de $\g(k)$ dans le groupe $M_k/M_{k-1}$ vaut $n_k$ pour tout $k$.
\end{lemm}
Enfin, pour tout $i$ compris entre 1 et $m$ notons $k_i$ le plus petit entier qui v\'erifie $a_{i}(k_i)\neq 0$. Quitte \`a effectuer une permutation des variables $X_i$, nous pouvons supposer que nous avons
$$1=k_1=k_2=\cdots=k_{i_0}<k_{i_0+1}\leq k_{i_0+2}\leq\cdots\leq k_m\leq g$$

\section[Arcs et d\'esingularisation plong\'ee torique]{Arcs et d\'esingularisation plong\'ee torique d'une hypersurface irr\'eductible \`a singularit\'e quasi-ordinaire}
Soit $\varphi(t):=(x_1(t),...,\,x_m(t),\,y(t))$ un arc trac\'e sur un germe $(X,\,0)$ d'hypersurface irr\'eductible \`a singularit\'e quasi-ordinaire d'exposants caract\'eristiques $a(1)$,..., $a(g)$ d\'efinie par un polyn\^ome de Weierstrass $f\in\C\{X_1,...,\,X_m\}[Y]$. Alors $f(x_1(t),...,\,x_m(t),\,y(t))=0$,
 donc $y(t)=\xi(x^{1/n}(t))$ o\`u $\xi=\xi_1+\cdots+\xi_g$ est une racine de $f$ (avec les notations du th\'eor\`eme \ref{quasi}) et les $x_i^{1/n}(t)$ sont des racines $n$-i\`emes de $x_i(t)$ dans $\C[[t^{1/n}]]$. Nous noterons souvent $\xi$ au lieu de $\xi(x^{1/n}(t))$ quand les racines $n$-i\`emes de $x_i(t)$ seront fix\'ees.\\
\begin{rema}\label{poids}Soit $X^a=X_1^{a_1}...X_m^{a_m}$ un mon\^ome de $\C[X_1^{1/n},...,\,X_m^{1/n}]$. Consid\'erons $m$ s\'eries $x_i(t)=\sum_kx_{i,k}t^k$ de $\C[[t]]$, et notons $l_i=\ord(x_i(t))$. Le choix d'une racine $n$-i\`eme de $x_i(t)$ dans $\C[[t^{1/n}]]$ d\'epend uniquement du choix d'une racine $n$-i\`eme de $x_{i,l_i}$ dans $\C$. Pour tout $i$ fixons une racine $n$-i\`eme de $x_{i,l_i}$ et notons la  $x_{i,l_i}^{1/n}$. Nous noterons alors sans \'equivoque  $x_{i,l_i}^{a_i}=(x_{i,l_i}^{1/n})^{na_i}$. Dans ce cas nous avons
$$x_i(t)^{a_i}=x_{i,l_i}^{a_i}t^{a_il_i}\left(1+\sum_{k\geq l_i+1}\frac{x_{i,k}}{x_{i,l_i}}t^{k-l_i}\right)^{a_1}=x_{i,l_i}^{a_i}t^{a_il_i}\left(1+\sum_{k\geq1}\frac{x_{i,k+l_i}}{x_{i,l_i}}t^{k}\right)^{a_1}$$
$$=x_{i,l_i}^{a_i}t^{a_il_i}\left(1+a_i\sum_{k\geq 1}\frac{x_{i,k+l_i}}{x_{i,l_i}}t^{k}+\cdots+\left(^{a_i}_j\right)\left(\sum_{k\geq 1}\frac{x_{i,k+l_i}}{x_{i,l_i}}t^{k}\right)^j+\cdots\right)$$
avec $\left(^{a_i}_j\right):=\frac{a_i(a_i-1)...(a_i-j+1)}{j!}$. Nous voyons que $x_1^{a_1}(t)...x_m^{a_m}(t)$ est dans $\C[[t]]$ si et seulement si $\sum_ia_il_i\in \N$. Si tel est le cas, pour tout $c\in \N$ avec $c>\sum_ia_il_i$, le coefficient de $t^c$ dans l'expression de la s\'erie $x_1^{a_1}(t)...x_m^{a_m}(t)$ est un polyn\^ome de la forme suivante :
$$\prod_{i=1}^mx_{i,l_i}^{a_i}P\left(x_{i,k}/x_{i,l_i};\,l_i+1\leq k\leq c-\sum_ja_jl_j+l_i\text{ et }1\leq i\leq m\right)$$
o\`u $P$ est un polyn\^ome quasi-homog\`ene de poids $c-\sum_ja_jl_j$ et o\`u $x_{i,k}/x_{i,l_i}$ est de poids $k-l_i$.

\end{rema}
Notons $b_k\, :\, \Z^m\lgw \Q$ la forme lin\'eaire
$$b_k(\underline{l}):=\sum_{i=1}^m a_i(k)l_i,\ \forall k\in\{0,...,\,g\}$$ et $M$ l'application lin\'eaire
$$\begin{array}{ccc} M\ :\ \Z^m\ & \lgw & \left(\Z/n\Z\right)^g\\
(l_1,...,\,l_m)& \lgm & (n\sum a_i(1)l_i,...,\, n\sum a_i(g)l_i)\end{array}$$
\\
Si $l_i$ est l'ordre de $x_i(t)$, alors d'apr\`es le th\'eor\`eme \ref{quasi}, nous voyons que n\'ecessairement $b_1(\underline{l})\in\N$ et donc $\xi_1(x^{1/n}(t))\in\C[[t]]$. En retranchant $\xi_1(x^{1/n}(t))$ \`a $y(t)$, nous voyons alors que $b_2(\underline{l})\in\N$ et donc que $\xi_2(x^{1/n}(t))\in\C[[t]]$. Par induction nous voyons que 
$$(\ord(x_1(t)),...,\,\ord(x_m(t)))\in\Ker\,M.$$
Inversement, si l'on se fixe $m$ s\'eries formelles en $t$, not\'ees $x_i(t)$ pour $1\leq i\leq m$, qui v\'erifient $(\ord(x_1(t)),...,\,\ord(x_m(t)))\in\Ker\,M\cap(\N^*)^m$, alors pour toute solution $\xi$ de $f$, nous avons $\xi(x^{1/n}(t))\in \C[[t]]$, et $(x_1(t),...,\,x_m(t),\,\xi(x^{1/n}(t)))$ d\'efinit un arc trac\'e sur $(X,\,0)$.\\
Nous allons maintenant relier $\Ker\,M$ aux r\'eseaux apparaissant dans une r\'esolution plong\'ee torique de $(X,\,0)$.
Nous rappelons ici la construction de la r\'esolution plong\'ee torique construite par P. Gonzalez Perez \cite{G-P2}. Dans la suite, nous noterons $Z_{\Sigma}$ la vari\'et\'e torique d'\'eventail $\Sigma$.\\
Soit $R=\C\{X_1,...,\,X_m\}[Y]/(f)$ l'anneau des fonctions de $(X,\,0)$. L'anneau $R$ est une $\C\{X_1,...,\,X_m\}$-alg\`ebre de finie. 
Nous avons alors un morphisme fini $\pi\,:\xymatrix{(X,\,0) \ar@{->>}[r] & (\A^{m}_{\C},\,0)}$ dont le lieu de ramification est inclus dans le germe \`a croisements normaux d\'efini par $X_1...X_m=0$ (non-ramifi\'e au-dessus du tore). Notons $\sigma$ l'\'eventail de $\A_{\C}^m$ \'gal \`a $\R_{\geq0}^m$ et $\D$ l'\'eventail $\sigma\oplus\R^g_{\geq0}=\R_{\geq0}^{m+d}\subset(N_{\D})_{\R}$ o\`u $N_{\D}$ est le r\'eseau $N\oplus\Z^g$. Soit $u_1,...,\,u_g$ la base canonique de $\{0\}\oplus\Z^g$. Les \'el\'ements de $\check{\D}\cap M_{\D}$ sont de la forme $(v,\,w)$ o\`u $v\in\check{\sigma}\cap M$ et $w=w_1u_1^*+\cdots+w_gu_g^*$ avec $w_i\in\Z_{\geq0}$. Nous allons consid\'erer le plongement $(X,\,0)\subset (Z_{\D},\,o_{\D})$, o\`u $o_{\D}$ est l'orbite de dimension 0 de la vari\'et\'e torique $Z_{\D}$, d\'efini par le morphisme de $\C\{X_1,...,\,X_m\}$-alg\`ebres :
$$\Psi : \C\{X_1,...,\,X_m\}[U_1,...,\,U_{g-1}]\lgw R$$
qui \`a $U_j$ associe $q_{j-1}(\xi)$ o\`u $(q_0,\,q_1,...,\,q_{g-1})$ est un syst\`eme de racines approch\'ees de $f$ (c'est-\`a-dire que $q_j$ est polyn\^ome quasi-ordinaire minimal de $\C\{X_1,...,\,X_m\}[Y]$ associ\'e \`a la branche $\xi_1+\cdots+\xi_j$ ; en particulier $q_0=Y$). Nous pouvons remarquer que $\Psi(U_1)=\xi$.\\
Soit 
$$\begin{array}{ccccc}
\varphi &:& M_{\D}& \lgw & M_g\\
 & & (v,\,w)& \lgm & v+\sum_{k=1}^gw_k\g(k)
\end{array}$$
Soit $L\subset(N_{\D})_{\R}$ l'espace lin\'eaire d\'efini comme \'etant l'orthogonal de $\Ker(\varphi)$ et soit $\Xi$ le c\^one $\D\cap L$. Nous appellerons \textit{squelette} de la singularit\'e le c\^one $\Xi$  (selon la d\'efinition de \cite{L-J-R1}).\\
Donnons tout d'abord la d\'efinition de r\'esolution partielle :
\begin{defi}
Soit $X$ une vari\'et\'e plong\'ee dans $Z_{\Sigma_1}$ une vari\'et\'e torique d'\'eventail $\Sigma_1$. Soit $\Sigma_2$ une subdivision de $\Sigma_1$. On dit que le morphisme torique $\varphi_{2,1} : Z_{\Sigma_2}\lgw Z_{\Sigma_1}$ est une r\'esolution plong\'ee partielle de $X$ si pour toute division r\'eguli\`ere $\Sigma_3$ de $\Sigma_2$ contenant tous les  c\^ones r\'eguliers de $\Sigma_2$, le morphisme torique induit $\varphi_{3,2} : Z_{\Sigma_3}\lgw Z_{\Sigma_2}$ compos\'e avec $\varphi_{2,1}$ est une r\'esolution plong\'ee de $X$.
\end{defi}
Nous avons alors le 
\begin{theo}\cite{G-P2}
Soit $\Sigma$ une division de $\D$ contenant $\Xi$. Alors le morphisme $\pi_{\Sigma}:Z_{\Sigma}\lgw Z_{\D}$ est une r\'esolution plong\'ee partielle du germe $(X,\,0)\subset (Z_{\D},\,o_{\D})$. \\
\end{theo}
Nous pouvons alors donner une description de $\Ker\,M$. Tout d'abord nous appelerons \textit{vecteur caract\'eristique} de $h$ (o\`u $h$ est un arc de $(\A^{m}_{\C},\,0)$ d\'efini par $(x_1(t),...,\,x_m(t)$) le vecteur $(x_1(t),...,\,x_m(t),\,q_0(\xi),...,\,q_g(\xi))$. 
\begin{prop}
Le squelette $\Xi$ correspond exactement \`a l'ensemble des vecteurs caract\'eristiques. En particulier le sous-r\'eseau $\Ker\,M\subset N$ est inclu dans la projection de $\Xi$ sur $N$.\\
\end{prop}

\begin{proof}
L'application lin\'eaire $\varphi$ est d\'efinie par la matrice 
$$\left(\begin{array}{cccccccc}
1 & 0 & 0& \cdots &0& \g_1(1)&\cdots&\g_1(g)\\
0& 1& 0& \cdots & 0& \vdots& \vdots&\vdots\\
\vdots & \vdots & \ddots & \vdots & \vdots& \vdots& \vdots&\vdots\\
0&0&\cdots &1&0& \vdots& \vdots&\vdots\\
0&0&\cdots&0&1& \g_m(1)&\cdots&\g_m(g)\\
\end{array}\right)$$
Le noyau de $\varphi$ est donc l'image de l'application lin\'eaire de $N_g$ dans $N_{\D}$ d\'efinie par la transpos\'ee de cette matrice. Or $N_g=\Ker\,M$. Donc $\Xi$ est exactement l'ensemble des \'el\'ements de la forme $(\underline{l},\,\underline{p})$ avec $\underline{l}\in\Ker\,M$ et $p_k=\sum_i\g_i(k)l_i$ pour $1\leq k\leq g$. Or, pour tout vecteur caract\'eristique $(x_1(t),...,\,x_m(t),\,q_0(\xi),...,\,q_g(\xi))$, nous avons $\ord(q_k(\xi))=\sum_i\g_i(k)\ord(x_i(t))$. Donc $\Xi$ correspond exactement \`a l'ensemble des vecteurs caract\'eristiques.\end{proof}

\section{D\'efinitions}
Nous allons \'etudier ici les s\'eries de Poincar\'e g\'eom\'etrique et arithm\'etique :
$$P_{geom,X(a(1),...,\,a(g)),0}(T):=\sum_{p\geq 0}[\pi_p(X(a(1),...,\,a(g))_{\infty})]T^p,$$
$$P_{arit,X(a(1),...,\,a(g)),0}(T):=\sum_{p\geq 0}\Chi_c(\varphi_p)T^p$$
o\`u $(X(a(1),...,\,a(g)),\,0)$ est un germe d'hypersurface irr\'eductible, mais non n\'ecessairement r\'eduit, d\'efini par  $f$  quasi-ordinaire. L'espace des arcs trac\'es sur le germe d'hypersurface d\'efini par $f$ et celui des arcs trac\'es sur le germe d'hypersurface r\'eduit associ\'e sont les m\^emes. Les espaces d'arcs tronqu\'es sont donc aussi les m\^emes. Nous pouvons donc supposer que le germe consid\'er\'e est r\'eduit. Supposons que $f$ est un polyn\^ome de Weierstrass de $\C\{X_1,...,\,X_m\}[Y]$ de degr\'e $n$.\\
Notons 
$$\Chi_{p,\,l_1,...,\,l_m}:=\pi_p(X(a(1),...,\,a(g))_{\infty})\cap\{(x_1(t),...,\,x_m(t),\,y(t))\,/\, \ord(x_i)=l_i\}$$
et 
$$\varphi_{p,\,l_1,...,\,l_m}:=\varphi_p\wedge (\ord(x_i)=l_i)$$
 avec $\varphi_p=\varphi_p(a(1),...,\,a(m))$ la formule
$$\exists z_1(t),...,\,z_m(t)\,\left((\ord(x_i-z_i)>p)\wedge (\ord(y-h(z))>p)\right)\ .$$
Nous noterons $X(a(1),...,\,a(g))_{p,cn}$ le sous-ensemble constructible de l'ensemble constructible $\pi_p(X(a(1),...,\,a(g))_{\infty})$ qui correspond \`a la troncation d'arcs pour lesquels un des $x_i(t)$ est nul et nous noterons $X(a(1),...,\,a(g))_{p,to}$ le sous-ensemble de $\pi_p(X(a(1),...,\,a(g))_{\infty})$ qui correspond \`a la troncation d'arcs qui ne peuvent pas se relever en arcs  pour lesquels un des $x_i(t)$ est nul. De m\^eme nous noterons $\varphi_{p,cn}(a(1),...,\,a(m))$ la formule 
$$\exists z_1(t),...,\,z_m(t)\,\left((\ord(x_i-z_i)>p)\wedge (\ord(y-h(z))>p)\wedge (\exists i\,/\, z_i(t)=0)\right)\ $$
 et $\varphi_{p,to}(a(1),...,\,a(m))$ la formule 
$$\big(\exists z_1(t),...,\,z_m(t)\,\left((\ord(x_i-z_i)>p)\wedge (\ord(y-h(z))>p)\right)\big)$$
$$\wedge\big(\neg\exists z_1(t),...,\,z_m(t)\,\left((\ord(x_i-z_i)>p)\wedge (\ord(y-h(z))>p)\wedge (\exists i\,/\, z_i(t)=0)\right)\big)\ $$ Nous avons \'evidemment
$$[\pi_p(X(a(1),...,\,a(g))_{\infty})]=[X(a(1),...,\,a(g))_{p,cn}]+[X(a(1),...,\,a(g))_{p,to}]$$
$$\text{ et }\Chi_c(\varphi_p(a(1),...,\,a(m)))=\Chi_c(\varphi_{p,to}(a(1),...,\,a(m)))+\Chi_c(\varphi_{p,cn}(a(1),...,\,a(m)))\ .$$
Notons aussi 
$$P_{geom}^{tore^c}(a(1),...,\,a(m))(T):=\sum_{p\geq0}[X(a(1),...,\,a(g))_{p,cn}]T^p\ ,$$
$$ P_{geom}^{tore}(a(1),...,\,a(m))(T):=\sum_{p\geq0}[X(a(1),...,\,a(g))_{p,to}]T^p\ ,$$
$$P_{arit}^{tore^c}(a(1),...,\,a(m))(T):=\sum_{p\geq0}\Chi_c(\varphi_{p,cn}(a(1),...,\,a(m)))T^p\ ,$$
$$\text{ et } P_{arit}^{tore}(a(1),...,\,a(m))(T):=\sum_{p\geq0}\Chi_c(\varphi_{p,to}(a(1),...,\,a(m)))T^p\ .$$
\begin{rema}
Si pour tout $i$ nous avons $l_i<+\infty$, alors $[\Chi_{p,\,l_1,...,\,l_m}]=0$ si et seulement si $\underline{l}\notin \Ker\,M$ (et de m\^eme $\Chi_c(\varphi_{p,\,l_1,...,\,l_m})=0$ si et seulement si $\underline{l}\notin \Ker\,M$).
\end{rema}

\section{\'Etude des troncations d'arcs ne se relevant pas en arcs pour lesquels un des $x_i(t)$ est nul}
Le terme $[X(a(1),...,\,a(g))_{p,to}]$ correspond aux arcs tronqu\'es qui ne peuvent pas se relever en arcs  pour lesquels un des $x_i(t)$ est nul. Si $(x_1(t),...,\,x_m(t),\,y(t))$ est un arc dont la $p$-troncation ne peut pas se relever en un arc pour lequel un des $x_i(t)$ est nul, nous notons $l_i=\ord(x_i(t))$. Nous avons alors la proposition suivante :
\begin{prop}\label{arctr}
Soit $h(t):=(x_1(t),...,\,x_m(t),\,y(t))$ un arc de $(X,\,0)$. Alors la $p$-troncation de $h$ ne peut pas se relever en un arc  pour lequel un des $x_i(t)$ est nul, si et seulement si l'une des deux conditions suivantes est v\'erifi\'ee :\\
\begin{itemize}
\item[C1)] Soit $ l_i\leq p$  pour tout $i$.\\
\item[C2)] Soit il existe $i$ tel que $+\infty>l_i\geq p+1$, et\\
\begin{itemize}
\item  $p-l_j\geq b_{k_i}(\underline{l})-b_{k_j}(\underline{l})$ pour tout $j$ tel que $k_j< k_i$,\\
\item  $b_{k_i}(\underline{l})\leq p$.
\end{itemize}
\end{itemize}
\end{prop}
\begin{rema}\label{bonplan}
En particulier, si $k_i\neq k_j$, alors $l_i$ et $l_j$ ne peuvent \^etre sup\'erieurs strictement \`a $p$ en m\^eme temps. En effet, supposons $l_i>p$ et $l_j>p$ avec $i\neq j$ et $k_i> k_j$. Alors nous avons selon la seconde condition : $p-l_j\geq b_{k_i}(\underline{l})-b_{k_j}(\underline{l})>0$, ce qui est contradictoire. Si il existe $i$ tel que $l_i>p$, nous noterons $r_{\underline{l}}$ l'entier $k_i$.
\end{rema}
\begin{proof} Montrons tout d'abord qu'un arc tronqu\'e qui ne peut pas se relever en arc  pour lesquel un des $x_i(t)$ est nul v\'erifie n\'ecessairement l'une des deux conditions. Tout d'abord, si $l_i\geq p+1$ et si $b_{k_i}(\underline{l})\geq p+1$, alors on peut relever cet arc tronqu\'e en un arc pour lequel $x_i(t)=0$, en consid\'erant un relev\'e quelconque $h$ de cet arc tronqu\'e, et en choisissant l'arc $h'$ dont toutes les coordonn\'ees, sauf la $i$-i\`eme, sont \'egales \`a celles de $h$, et la $i$-\`eme est nulle.\\
Maintenant supposons que  $l_i\geq p+1$ et que $p-l_j< b_{k_i}(\underline{l})-b_{k_j}(\underline{l})$ pour un $j$ tel que $k_j< k_i$. En particulier, d'apr\`es la remarque \ref{bonplan}, pour $j$ tel que $k_j\neq k_i$, $l_j\leq p$ et donc les $x_{j,k}$ sont fix\'es pour $k\leq p$ car les $x_j(t)$ sont fix\'es modulo $t^{p+1}$. Consid\'erons alors le coefficient de $t^{b_{k_i}(\underline{l})}$ dans l'\'ecriture de $\xi_{k_j}(x^{1/n})$ et dans l'\'ecriture de $\xi_{k_i}(x^{1/n})$. Dans  $\xi_{k_j}(x^{1/n})$, celui-ci est de la forme 
$$\prod_{r\neq j}x_{r,l_r}^{a_{r}(k_j)}.\left(\frac{x_{j,b_{k_i}(\underline{l})-b_{k_j}(\underline{l})+l_j}}{x_{j,l_j}}+P\left(\frac{x_{r,k}}{x_{r,l_r}};\, k\leq b_{k_i}(\underline{l})-b_{k_j}(\underline{l})+l_r\text{ et }1\leq r\leq m\right)\right)$$
o\`u $P$ ne d\'epend pas de $\frac{x_{j,b_{k_i}(\underline{l})-b_{k_j}(\underline{l})+l_j}}{x_{j,l_j}}$. Comme $b_{k_i}(\underline{l})-b_{k_j}(\underline{l})+l_j>p$, le terme $\frac{x_{j,\,b_{k_i}(\underline{l})-b_{k_j}(\underline{l})+l_j}}{x_{j,l_j}}$ n'apparait pas dans l'\'ecriture de $x_j(t)$ modulo $t^{p+1}$. Il n'apparait pas non plus dans l'expression des coefficients de $t^d$, pour $d<b_{k_i}(\underline{l})$, dans l'\'ecriture de $\xi(x^{1/n})$. 
D'autre part, le coefficient de $t^{b_{k_i}(\underline{l})}$ dans l'\'ecriture de $\xi_{k_i}(x^{1/n})$ est de la forme :
$$\prod_{r\neq i}x_{r,l_r}^{a_{r}(k_i)}.x_{i,l_i}^{a_{r}(k_i)}$$
En particulier, nous pouvons choisir $x_{i,l_i}$ \'egal \`a z\'ero en donnant la bonne valeur au terme $\frac{x_{j,b_{k_i}(\underline{l})-b_{k_j}(\underline{l})+l_j}}{x_{j,l_j}}$. Ceci n'affecte alors ni la valeur de $x_j(t)$ modulo $t^{p+1}$, ni la valeur des coefficients de $t^d$, pour $d<b_{k_i}(\underline{l})$, dans l'\'ecriture de $\xi(x^{1/n})$. Nous pouvons continuer ainsi par r\'ecurrence croissante sur $c$ en regardant le coefficient de $t^c$ dans l'\'ecriture de $\xi(x^{1/n})$ et annuler le terme $x_{i,c-\sum_{r\neq i}a_r(k_i)l_r}$ en modifiant \'eventuellement le terme $\frac{x_{j,c-b_{k_j}(\underline{l})+l_j}}{x_{j,l_j}}$.\\
\\
Montrons maintenant la suffisance de ces deux conditions. La premi\`ere condition implique clairement que l'arc tronqu\'e ne peut pas se relever en un arc dont l'une des coordonn\'ees est nulle. Montrons que la seconde est aussi suffisante.\\
Soit $h$ un arc qui v\'erifie la seconde condition. Tout d'abord, comme $l_j\leq p$ pour $j$ tel que $k_j\neq k_i$, les coordonn\'ees $x_j(t)$ d'un rel\`evement d'une $p$-troncation de $h$ sont non nulles. Plus particuli\`erement, pour $j$ tel que $k_j\neq k_i$, $x_{j,k}$ est fix\'e pour $l_j\leq k\leq p$. Le coefficient de $t^{b_{1}(\underline{l})}$ fixe $\prod_j x_{j,l_j}^{a_j(1)}$. Par r\'ecurrence croissante sur les coefficients de $t^c$, pour $c<b_{k_i}(\underline{l})$, $\prod_j x_{j,l_j}^{a_j(r)}$ est fix\'e pour $r<k_i$ car seuls les $x_{j,k}$ pour $j\neq i$ et $k\leq c-b_{k_j}(\underline{l})+l_j<p$ apparaissent dans son expression.
Consid\'erons alors le coefficient de $t^{b_{k_i}(\underline{l})}$ dans l'expression  de la coordonn\'ee $y(t)$ de $h$. Nous avons $y(t)=\xi_1(x^{1/n}(t))+\cdots+\xi_{k_i}(x^{1/n}(t))+\cdots+\xi_g(x^{1/n}(t))$ pour un choix d'une racine $n$-i\`eme des $x_j(t)$. Le coefficient de $t^{b_{k_i}(\underline{l})}$ dans l'expression de $\xi_r(x^{1/n}(t))$ est de la forme (pour $r<k_i$) :
$$\prod_{j\neq i}x_{j,l_j}^{a_j(r)}P_r(x_{j,k}/x_{j,l_j}, k\leq b_{k_i}(\underline{l})-b_{k_r}(\underline{l})+l_j)$$ d'apr\`es la remarque \ref{poids}. Or nous avons $p-l_j\geq b_{k_i}(\underline{l})-b_{k_j}(\underline{l})\geq b_{r}(\underline{l})-b_{k_j}(\underline{l})$, et donc $P_r(x_{j,k}/x_{j,l_j}, k\leq b_{k_i}(\underline{l})-b_{k_r}(\underline{l})+l_j)$ est fix\'e du fait que les $x_j(t)$ sont fix\'es modulo $t^{p+1}$. D'apr\`es la r\'ecurrence, le coefficient consid\'er\'e est donc fix\'e. Le coefficient de $t^{b_{k_i}(\underline{l})}$ dans l'expression de $\xi_{k_i}(x^{1/n}(t))$ est de la forme :
$$\prod_{j\neq i}x_{j,l_j}^{a_j(k_i)}.x_{i,l_i}^{a_i(k_i)}$$
Donc le coefficient de $t^{b_{k_i}(\underline{l})}$ dans l'expression de $y(t)$ est de la forme $a x_{i,l_i}^{a_i(k_i)}+b$ o\`u $b$ est fix\'e et $a$ est non nul par hypoth\`ese sur $h$. Comme $x_{i,l_i}\neq 0$ par hypoth\`ese sur $h$, $b\neq 0$ et n\'ecessairement tout relev\'e de la $p$-troncation de $h$ a sa $i$-i\`eme coordonn\'ee non nulle. Plus pr\'ecis\'ement, nous pouvons remarquer  que n\'ecessairement $b_{k_i}(\underline{l})=\sum_ja_j(k_i)l_j$ est fix\'e.\\
\end{proof}
Nous avons donc 
\begin{equation}\label{union}X(a(1),...,\,a(g))_{p,to}=\bigcup_{(l_1,...,\,l_m)\in D(m)_p}\Chi_{p,\,l_1,...,\,l_m}\end{equation} o\`u 
$$\Sigma_{m,p}:=\left\{(l_1,...,\,l_m)\,/\, l_i\leq p\right\}\bigcup$$
$$\left(\bigcup_{i=1}^m\left\{(l_1,...,\,l_m)\,/\,l_i>p,\,p-l_j\geq b_{k_i}(\underline{l})-b_{k_j}(\underline{l}),\,k_j< k_i,\,b_{k_i}(\underline{l})\leq p\right\}\right)\,$$
$$\text{et } D(m)_p:=\Sigma_{m,p}\cap\Ker\,M\cap (\N^*)^m.$$
Il nous suffit donc de calculer $\Chi_{p,\,l_1,...,\,l_m}$ et $\varphi_{p,\,l_1,...,\,l_m}$, et d'\'etudier l'union (\ref{union}).

\section{Calcul de $\Chi_{p,\,l_1,...,\,l_m}$ et $\varphi_{p,\,l_1,...,\,l_m}$}Pour calculer $\Chi_{p,\,l_1,...,\,l_m}$ et $\varphi_{p,\,l_1,...,\,l_m}$, nous allons voir que $\Chi_{p,\,l_1,...,\,l_m}$ peut s'\'ecrire sous la forme $\ovl{\Chi}_{p,\,l_1,...,\,l_m}\times \A_{\C}^{n(p,\underline{l})}$, o\`u $n(p,\underline{l})$ est un entier, et nous allons trouver un rev\^etement galoisien $W\lgw \ovl{\Chi}_{p,\,l_1,...,\,l_m}$ surjectif o\`u $[W]$ est facile \`a calculer et l'action du groupe de Galois du rev\^etement est facile \`a d\'ecrire.\\
Pour cela, notons $x_j=\sum_{r\geq l_j}x_{j,r}t^r$ et $z_j=z_{j,0}t^{\frac{l_j}{n}}(1+\sum_{r\geq 1}z_{j,r}t^r)$. Les termes $z_{j,r}$ pour $r>p-l_j$ n'apparaissent pas dans l'expression de $z_j^n$ modulo $t^{p+1}$ mais $z_{j,p-l_j}$ y apparait. De m\^eme les termes $z_{j,r}$ pour $r>p-b_{k_j}(\underline{l})$ n'apparaissent pas dans l'expression de $\xi(z)$ modulo $t^{p+1}$ si $p\geq b_{k_j}(\underline{l})$ mais $z_{j,p-b_{k_j}(\underline{l})}$ y apparait. Si  $p< b_{k_j}(\underline{l})$, alors aucun terme $z_{j,r}$ n'apparait dans l'expression de $\xi(z)$ modulo $t^{p+1}$.\\
Soit $$W_{\underline{l}}:=\mathbb{G}_{m,\C}^m\times \A_{\C}^{\max\{p-l_1,\,p-b_{k_1}(\underline{l})\}}\times\cdots\times\A_{\C}^{\max\{p-l_m,\,p-b_{k_m}(\underline{l})\}}.$$
On consid\`ere sur $W_{\underline{l}}$ les coordonn\'ees $(z_{1,0},...,\,z_{m,0})$ sur le premier facteur et les coordonn\'ees $(z_{j,1},...,\,z_{j,\max\{p-l_j,\,p-b_{k_j}(\underline{l})\}})$ sur le $(j+1)$-i\`eme facteur.\\
Soit $V$ la vari\'et\'e $\A^{p(m+1)}$ ; consid\'erons le morphisme de $\C$-sch\'emas $h_{\underline{l}}\ :\ W_{\underline{l}}\lgm V$ qui envoie le point de coordonn\'ees pr\'ec\'edentes sur les $p$ premiers coefficients de $z_1^n$,..., $z_m^n$ et $\xi(z)$ (o\`u $z_j(t):=z_{j,0}t^{\frac{l_i}{n}}(1+\sum_{r=1}^{\max(p-l_j,\,p-b_{k_j}(\underline{l}))}z_{j,r}t^r)$) :
$$ h_{\underline{l}}\, :\, W_{\underline{l}}\longrightarrow V$$
$$ z_i=z_{i,0}t^{\frac{l_i}{n}}(1+\sum_{j=1}^{\max(p-l_i,\,p-b_{k_i}(\underline{l}))}z_{i,j}t^j)\lgm z_1^n,...,\,z_m^n,\,\xi(z)\qquad\qquad\qquad\quad$$
L'image de $h_{\underline{l}}$ co\"incide clairement avec $\Chi_{p,\,l_1,...,\,l_m}$, mais ce morphisme n'est pas fini en g\'en\'eral.\\
Nous allons noter (si $I_{p,\,\underline{l}}$ est l'ensemble des indices $j$ pour lesquels $l_j>p$)
$$D_{0,i}(m)_p:=\left\{\underline{l}\in D(m)_p\ \backslash\ l_j\leq p \ \ \forall j,\ l_i-b_{k_i}(\underline{l})\geq l_j-b_{k_j}(\underline{l}) \ \ \forall j\neq i \right.$$
$$\left.\text{ et }l_i-b_{k_i}(\underline{l})> l_j-b_{k_j}(\underline{l})\  \ \forall j<i\right\}$$
et
$$D_{q,i}(m)_p:=\left\{\underline{l}\in D(m)_p\ \backslash\ l_i> p, \ 
Card(I_{p,\,\underline{l}})=q,\ l_i-b_{k_i}(\underline{l})\geq l_j-b_{k_j}(\underline{l}) \ \ \forall j\neq i \right.$$
$$\left.\text{ et }l_i-b_{k_i}(\underline{l})> l_j-b_{k_j}(\underline{l})\  \ \forall j<i\right\}.$$
\begin{rema}\label{bonplan2}Si $l_i>p$ et si $j$ est tel que $k_j\neq k_i$, alors $l_i-b_{k_i}(\underline{l})>l_j-b_{k_j}(\underline{l})$. En effet, si $k_i>k_j$, alors par hypoth\`ese nous avons $p-l_j\geq b_{k_i}(\underline{l})-b_{k_j}(\underline{l})$. Or $p-l_j<l_i-l_j$ et le r\'esultat s'ensuit. Si $k_i<k_j$, alors $0<l_i-p\leq l_i-l_j$. Or $b_{k_i}(\underline{l})-b_{k_j}(\underline{l})<0$ et le r\'esultat s'ensuit l\`a encore. En particulier nous avons
$$D(m)_p=\coprod_{0\leq i,\,q\leq m}D_{i,q}(m)_p.$$
\end{rema}
\subsection{Cas C1}
Soit $\underline{l}\in D_{0,i}(m)_p$. Nous notons alors
$$W'_{\underline{l}}:=\mathbb{G}_{m,\C}^m\times \A_{\C}^{p-l_1}\times\cdots\times\A_{\C}^{\max\{p-l_i,\,p-b_{k_i}(\underline{l})\}}\times\cdots\times\A_{\C}^{p-l_m}$$
et
$$W''_{\underline{l}}:=\mathbb{G}_{m,\C}^m\times \A_{\C}^{p-l_1}\times\cdots\times\A_{\C}^{p-l_m}.$$
On consid\`ere sur $W'_{\underline{l}}$ les coordonn\'ees $(z_{1,0},...,\,z_{m,0})$ sur le premier facteur, les coordonn\'ees $(z_{j,1},...,\,z_{j,p-l_j})$ sur le $(j+1)$-i\`eme facteur pour $j\neq i$ et les coordonn\'ees $(z_{i,1},...,\,z_{i,\max\{p-l_i,\,p-b_{k_i}(\underline{l})\}})$ sur $(i+1)$-i\`eme facteur. On consid\`ere sur $W''_{\underline{l}}$ les coordonn\'ees $(z_{1,0},...,\,z_{m,0})$ sur le premier facteur, et les coordonn\'ees $(z_{j,1},...,\,z_{j,p-l_j})$ sur le $(j+1)$-i\`eme facteur pour tout $j$.\\
Nous avons clairement $$W''_{\underline{l}}\subset W'_{\underline{l}}\subset W_{\underline{l}}.$$
Nous notons 
$$e:=\max_j\{l_j-b_{k_j}(\underline{l}),\ 0\}=\max\{l_i-b_{k_i}(\underline{l}),\ 0\}.$$
Par ailleurs, nous notons $k_{\underline{l}}$ l'unique entier qui v\'erifie les in\'egalit\'es suivantes : 
\begin{equation}b_{k_{\underline{l}}}(\underline{l})=\sum_j a_j(k_{\underline{l}})l_j\leq p-e<\sum_j a_{j}(k_{\underline{l}}+1)l_j=b_{k_{\underline{l}}+1}(\underline{l}).\end{equation}
Nous avons alors le r\'esultat suivant :
\begin{lemm}
 Le morphisme $h_{\underline{l}}$ restreint \`a $W'_{\underline{l}}$ est fini et d'image $\Chi_{p,l_1,...,l_m}$. Si nous notons $\ovl{\Chi}_{p,l_1,...,l_m}:=h_{\underline{l}}(W''_{\underline{l}})$, nous avons $\Chi_{p,\,l_1,...,\,l_m}=\ovl{\Chi}_{p,l_1,...,l_m}\times \A_{\C}^{l_i-b_{k_i}(\underline{l})}$ et $h_{\underline{l}}$ restreint \`a $W''_{\underline{l}}$ est un rev\^etement galoisien de $\ovl{\Chi}_{p,\,l_1,...,\,l_m}$. Son groupe de Galois est le groupe $G_{k_{\underline{l}}}$, sous-groupe commutatif de $\U_n^m$ form\'e des \'el\'ements $(\e_1,...,\,\e_m)$ qui v\'erifient les \'equations $\prod_i \e_i^{na_r(i)}=1$ pour $ r=1,..,k_{\underline{l}}$ o\`u $\U_n$ est le groupe des racines $n$-i\`emes de l'unit\'e, qui agit par multiplication terme \`a terme sur $\mathbb{G}_{m,\C}^m$.
\end{lemm}
\begin{proof}
Nous allons d'abord montrer que l'image de $W'_{\underline{l}}$ par $h_{\underline{l}}$ est \'egale \`a $\Chi_{p,l_1,...,l_m}$.\\
Soit $z=(z_1,...,\,z_m)\in W_{\underline{l}}$, et montrons que l'on peut trouver $w\in W'_{\underline{l}}$  ayant m\^eme image par $h_{\underline{l}}$ que $z$. N\'ecessairement nous devons avoir
$$z_{j,0}^n=w_{j,0}^n \text{  si }\ l_j\leq p,$$
$$\left(1+\sum_{k=1}^{p-l_j}z_{j,k}t^k\right)^n=\left(1+\sum_{k=1}^{p-l_j}w_{j,k}t^k\right)^n\text{ mod } t^{p-l_j+1}\ \text{si }\ l_j\leq p.$$\\
Comme $l_j\leq p$ pour tout $j$, n\'ecessairement $z_{j,0}^n=w_{j,0}^n$ et $z_{j,k}=w_{j,k}$, pour tout $j$ et $1\leq k\leq p-l_j$. Nous posons alors $w_{j,k}=z_{j,k}$, pour tout $j$ et $0\leq k\leq p-l_j$. Les coefficients de $t^c$ dans l'\'ecriture de $\xi(z)$, pour $c\leq p-e$, ne d\'ependent que des $z_{j,k}$ pour $k\leq p-l_j$.  Le coefficient de $t^{p-e+1}$ est de la forme 
$$\prod_jz_{j,l_j}^{a_j(k_j)}.z_{i, p-l_i+1}+P$$
o\`u $P$ ne d\'epend que des $z_{i,k}$ pour $k\leq p-l_i$ et des $z_{j,k}$ pour $j\neq i$ et $k\leq p-l_j+1$. Comme $\prod_jz_{j,l_j}^{a_j(k_j)}=\prod_jw_{j,l_j}^{a_j(k_j)}\neq 0$, nous pouvons trouver $w_{i, p-l_i+1}$ de telle mani\`ere \`a choisir (s'il apparaissent ici) les $w_{j, p-l_j+1}$ \'egaux \`a z\'ero pour $j\neq i$ sans changer la valeur de ce coefficient. Nous pouvons continuer ainsi par r\'ecurrence croissante sur les coefficients de $t^c$ dans l'\'ecriture de $\xi(z)$, pour voir que l'on peut trouver $w\in W'_{\underline{l}}$ tel que $h_{\underline{l}}(z)=h_{\underline{l}}(w)$, le syst\`eme d'\'equations apparaissant ici \'etant triangulaire. Le morphisme $h_{\underline{l}}$ restreint \`a $W'_{\underline{l}}$ est donc fini. \\
Nous remarquons au passage, que le syst\`eme d'\'equations obtenu des coefficients de $t^c$ dans l'\'ecriture de $\xi(z)$ est triangulaire en les variables $z_{i,k}$ pour $p-l_i+1\leq k\leq p-b_{k_i}(\underline{l})$, pour $p-e+1\leq c\leq p$. Donc nous pouvons \'ecrire $\Chi_{p,l_1,...,l_m}$ sous la forme $\ovl{\Chi}_{p,l_1,...,l_m}\times \A_{\C}^{l_i-b_{k_i}(\underline{l})}$ o\`u $\ovl{\Chi}_{p,l_1,...,l_m}=h_{\underline{l}}(W''_{\underline{l}})$.\\
D\'eterminons la fibre au-dessus d'un point de $\ovl{\Chi}_{p,l_1,...,l_m}$. Remarquons que $k_i\leq k_{\underline{l}}$ car $b_{k_{\underline{l}}+1}(\underline{l})> p-e=p-l_i+b_{k_i}(\underline{l})\geq  b_{k_i}(\underline{l})$. Soient $z$ et $w$ dans $W''_{\underline{l}}$ tels que $h_{\underline{l}}(z)=h_{\underline{l}}(w)$. Dans ce cas $z_{j,0}^n=w_{j,0}^n$ pour tout $j$ et $z_{j,k}=w_{j,k}$ pour tout $j$ et pour tout $1\leq k\leq p-l_j$. Consid\'erons alors le coefficient de $t^c$ dans l'\'ecriture de $\xi(z)=\xi(w)$. Pour $c=b_{k_1}(\underline{l})$, ce coefficient nous  permet de dire que $(z_{1,0},...,z_{m,0})=\e(w_{1,0},...,w_{m,0})$ pour un $\e\in G_{k_1}$. Pour $b_{k_1}(\underline{l})<c<b_{k_2}(\underline{l})$, ces coefficients n'apportent aucune information suppl\'ementaire. Pour $c=b_{k_2}(\underline{l})$, coefficient nous  permet de dire que $(z_{1,0},...,z_{m,0})=\e(w_{1,0},...,w_{m,0})$ pour un $\e\in G_{k_2}$.
Par r\'ecurrence sur $c$, nous voyons alors que $(z_{1,0},...,z_{m,0})=\e(w_{1,0},...,w_{m,0})$ pour un $\e\in G_{k_{\underline{l}}}$. Comme $l_j<+\infty$ pour tout $j$, le rev\^etement est \'etale et donc galoisien de groupe de Galois $G_{k_{\underline{l}}}$.\\\end{proof}
Nous avons alors 
$$[\Chi_{p,\,l_1,...,\,l_m}]=[\mathbb{G}_{m,\C}^m/G_{k_{\underline{l}}}]\DL^{pm+\max\{l_i-b_{k_i}(\underline{l}),\,0\}-\sum_{j=1}^ml_j}\ .$$

\subsection{Cas C2}
Soit $\underline{l}\in D_{q,i}(m)_p$ avec $q\geq 1$. Nous notons alors
$$W'_{\underline{l}}:=\mathbb{G}_{m,\C}^{m-q+1}\times\left(\prod_{j\notin I_{\underline{l},p}} \A_{\C}^{p-l_j}\right)\times\A_{\C}^{p-b_{k_i}(\underline{l})}$$
et
$$W''_{\underline{l}}:=\mathbb{G}_{m,\C}^{m-q+1}\times\left(\prod_{j\notin I_{\underline{l},p}} \A_{\C}^{p-l_j}\right).$$
On consid\`ere sur $W'_{\underline{l}}$ les coordonn\'ees $(z_{j_1,0},...,\,z_{j_{m-q},0})$ sur le premier facteur (o\`u $\{j_1,...,\,j_{m-q+1}\}=\{1,...,\,m\}\backslash \left\{I_{\underline{l},p}\backslash\{i\}\right\}$), les coordonn\'ees $(z_{j,1},...,\,z_{j,p-l_j})$ sur le $(j+1)$-i\`eme facteur pour $j\notin I_{\underline{l},p}$ et les coordonn\'ees $(z_{i,1},...,\,z_{i,p-b_{k_i}(\underline{l})})$ sur $(i+1)$-i\`eme facteur. On consid\`ere ces deux vari\'et\'es plong\'ees dans $W_{\underline{l}}$ en identifiant un \'el\'ement de coordonn\'ees 
$$(z_{j_1,0},...,\,z_{j_{m-q+1},0},\ z_{j,1},...,\,z_{j,p-l_j}, j\notin I_{\underline{l},p})$$
avec l'\'el\'ement de $W_{\underline{l}}$ de coordonn\'ees
$$(z_{1,0},...,\,z_{m,0},\ z_{j,1},...,\,z_{j,p-l_j}, j\neq i, z_{i,1},...,\,z_{i,p-b_{k_i}(\underline{l})})$$
en posant $z_{j,0}=1$  et $z_{j,k}=0$  si $j\in I_{p,\,\underline{l}}$ et $k\leq p-l_j$.
 Nous avons alors le r\'esultat suivant :
\begin{lemm}
Le morphisme $h_{\underline{l}}$ restreint \`a $W'_{\underline{l}}$ est fini et d'image $\Chi_{p,\,l_1,...,\,l_m}$. Si nous notons $\ovl{\Chi}_{p,\,l_1,...,\,l_m}:=h_{\underline{l}}(W''_{\underline{l}})$, nous avons $\Chi_{p,\,l_1,...,\,l_m}=\ovl{\Chi}_{p,\,l_1,...,\,l_m}\times \A_{\C}^{l_i-b_{k_i}(\underline{l})}$ et $h_{\underline{l}}$ restreint \`a $W''_{\underline{l}}$ est un rev\^etement galoisien. Son groupe de Galois est le groupe $G_{k_i}\cap \U^{m-q+1}_n$, o\`u $\U_n^{m-q+1}$ est le sous-groupe de $\U_n^m$ dont les \'el\'ements ont les coordonn\'ees d'indice dans $I_{\underline{l},p}\backslash\{i\}$ \'egales \`a 0, et $G_{k_i}$ est le sous-groupe de $\U_n^m$ form\'e des \'el\'ements $(\e_1,...,\,\e_{m-q+1})$ qui v\'erifient les \'equations $\prod_{j\notin I_{p,\,\underline{l}}\backslash\{i\}} \e_j^{na_r(j)}=1$ pour $ r=1,..,k_{i}$ o\`u $\U_n$ est le groupe des racines $n$-i\`emes de l'unit\'e, qui agit par multiplication terme \`a terme sur $\mathbb{G}_{m,\C}^{m-q}$. 
\end{lemm}
\begin{proof}
Nous allons d'abord montrer que l'image de $W'_{\underline{l}}$ par $h_{\underline{l}}$ est \'egale \`a $\Chi_{p,l_1,...,l_m}$.\\
Soit $z=(z_1,...,\,z_m)\in W_{\underline{l}}$, et montrons que l'on peut trouver $w\in W'_{\underline{l}}$  ayant m\^eme image par $h_{\underline{l}}$ que $z$. En particulier nous avons
$$z_{j,0}^n \text{ est fix\'e  si }\ l_j\leq p,$$
$$\left(1+\sum_{k=1}^{p-l_j}z_{j,k}t^k\right)^n\text{ mod } t^{p-l_j+1}\ \text{est fix\'e   si }\ l_j\leq p.$$\\
Les coefficients de $t^c$ dans l'\'ecriture de $\xi(z)$, pour $c< b_{k_i}(\underline{l})$, ne d\'ependent que des $z_{j,k}$ pour $j\notin I_{p,\,\underline{l}}$ et $k\leq p$, d'apr\`es la remarque \ref{bonplan2}. Si $c=b_{k_i}(\underline{l})$, ce coefficient est de la forme
$$\prod_{j\neq i}z_{j,0}^{a_j(k_i)}.z_{i,0}^{a_i(k_i)}+P$$
o\`u $P$ ne d\'epend que des $z_{j,k}$ pour $j\notin I_{\underline{l},p}$, $k\leq p-l_j+1$. Si nous posons $w_{j,0}=1$ pour $j\in I_{p,\,\underline{l}}\backslash \{i\}$ et $w_{j,p-l_j+1}=0$, il existe toujours un $w_{i,0}$ qui permet de garder la valeur de ce coefficient inchang\'ee. Nous pouvons continuer ainsi par r\'ecurrence croissante sur les coefficients de $t^c$ dans l'\'ecriture de $\xi(z)$, pour voir que l'on peut trouver $w\in W'_{\underline{l}}$ tel que $h_{\underline{l}}(z)=h_{\underline{l}}(w)$, le syst\`eme d'\'equations apparaissant ici \'etant triangulaire. Le morphisme $h_{\underline{l}}$ restreint \`a $W'_{\underline{l}}$ est donc fini. \\
Nous remarquons au passage, que le syst\`eme d'\'equations obtenu des coefficients de $t^c$ dans l'\'ecriture de $\xi(z)$ est triangulaire en les variables $z_{i,k}$, pour $0\leq k\leq p-b_{k_i}(\underline{l})$ et $b_{k_i}(\underline{l})\leq c\leq p$. Donc nous pouvons \'ecrire $\Chi_{p,l_1,...,l_m}$ sous la forme $\ovl{\Chi}_{p,l_1,...,l_m}\times \A_{\C}^{l_i-b_{k_i}(\underline{l})}$ o\`u $\ovl{\Chi}_{p,l_1,...,l_m}=h_{\underline{l}}(W''_{\underline{l}})$.\\
Comme pr\'ec\'edemment nous voyons que deux \'el\'ements $z$ et $w$ de $W''_{\underline{l}}$ ont m\^eme image par $h_{\underline{l}}$ si et seulement si $z=\e w$ o\`u $\e\in G_{k_i}\cap \U^{m-q+1}_n$.\\
\end{proof}
Nous avons alors 

$$[\Chi_{p,\,l_1,...,\,l_m}]=[\mathbb{G}_{m,\C}^{m-q+1}/G'_{k_i}]\DL^{p(m-q+1)-\sum_{j\notin I_{p,\,\underline{l}}}l_j-b_{k_i}(\underline{l})}\ .$$
\subsection{Calcul de $\Chi_{p,\,l_1,...,\,l_m}$ et $\varphi_{p,\,l_1,...,\,l_m}$}
Nous allons maintenant \'enoncer un lemme utile pour achever le calcul :
\begin{lemm}
Soit $G$ un sous-groupe de $\prod_{i=1}^m\U_{n_i}$ agissant sur $\mathbb{G}_{m,\C}^m$ par multiplication sur chaque terme. Alors la vari\'et\'e quotient $\mathbb{G}_{m,\C}^m/G$ est isomorphe \`a $\mathbb{G}_{m,\C}^m$.
\end{lemm}
\begin{proof}
D'un point de vue torique, le c\^one associ\'e \`a $\mathbb{G}_{m,\C}^m$ est l'origine de $\Z^m$. Le semi-groupe de cette vari\'et\'e est donc $\Z^m$ en entier. L'ensemble des puissances des mon\^omes invariants par l'action de $G$ forment un semi-groupe $N$ de $\Z^m$ qui est en fait un groupe car $X_1^{k_1}...X_m^{k_m}$ est invariant sous l'action de $G$ si et seulement si $X_1^{-k_1}...X_m^{-k_m}$ l'est aussi. Donc $N$ est isomorphe \`a $\Z^l$ pour $l\leq m$ et $N\otimes_{\Z}\R\equiv \R^l$. Or pour tout $i$, $X_i^{n_i}$ est invariant par $G$ donc $N\otimes_{\Z}\R=\R^m$ et $l=m$. Donc  $\mathbb{G}_{m,\C}^m/G$ est isomorphe \`a $\mathbb{G}_{m,\C}^m$.\\
\end{proof}

Nous avons donc $[\mathbb{G}_{m,\C}^r/G_{k}]=(\DL-1)^r$ et nous obtenons :
$$\text{pour le cas C1 : }\ [\Chi_{p,\,l_1,...,\,l_m}]=(\DL-1)^m\DL^{pm}\DL^{-\sum_{j=1}^ml_j+\max\{l_i-b_{k_i}(\underline{l})),\,0\}}.$$
$$\text{pour le cas C2 : }\ [\Chi_{p,\,l_1,...,\,l_m}]=(\DL-1)^{m-q+1}\DL^{p(m-q+1)-\sum_{j\notin I_{\underline{l},p}}l_j-b_{k_i}(\underline{l})}.$$
\\
Nous utilisons le m\^eme morphisme $h_{\underline{l}}$ pour calculer $\Chi_c(\varphi_{p,\,l_1,...,\,l_m})$. En effet, dans le cas C1, la mesure arithm\'etique de l'image de $h_{\underline{l}}$ co\"incide avec $\Chi_c(\varphi_{p,\,l_1,...,\,l_m})$. Or $\varphi_{p,\,l_1,...,\,l_m})=\ovl{\varphi}_{p,\,l_1,...,\,l_m}\wedge\varphi'$ o\`u $\ovl{\varphi}_{p,\,l_1,...,\,l_m}$ est la formule qui d\'efinit $\ovl{\Chi}_{p,\,l_1,...,l_m}$ et $\varphi'$ est la formule qui d\'efinit $\A_{\C}^{\max\{l_i-b_{k_i}(\underline{l}),\,0\}}$. La mesure arithm\'etique de $W''_{\underline{l}}$ est \'egale \`a $\DL^{p(m-q)-\sum_{j\notin I_{p,\,\underline{l}}}l_j}$ car la formule qui d\'efinit $W''_{\underline{l}}$ est sans quantificateurs. Le morphisme $h_{\underline{l}}$ restreint \`a $W''_{\underline{l}}$ est un rev\^etement galoisien de groupe de Galois fini, de cardinal $n/n_{k_{\underline{l}}}$ dans le cas C1 et de cardinal $n(i,\,I_{p,\underline{l}})$ dans le cas C2 ($n(i,\,I_{p,\underline{l}})$ ne d\'epend que de $i$ et de $I_{p,\underline{l}})$). La formule $\varphi'$ est sans quantificateurs, donc sa mesure vaut $\DL^{l_i-b_{k_i}(\underline{l})}$. Donc nous avons 
$$\text{cas C1 : }\ \Chi_c(\varphi_{p,\,l_1,...,\,l_m})=\frac{n_{k_{\underline{l}}}}{n}(\DL-1)^m\DL^{pm}\DL^{-\sum_{j=1}^ml_j+\max\{l_i-b_{k_i}(\underline{l})),\,0\}}.\qquad$$
$$\text{cas C2 : }\ \Chi_c(\varphi_{p,\,l_1,...,\,l_m})=\frac{1}{n(i,\,I_{p,\underline{l}})}(\DL-1)^{m-q+1}\DL^{p(m-q+1)-\sum_{j\notin I_{p,\,\underline{l}}}l_j-b_{k_i}(\underline{l})}.$$\\

\section{\'Etude de l'union (\ref{union})}
L'union (\ref{union}) n'est cependant pas toujours disjointe comme le montre l'exemple ci-dessous :
\begin{exem}
Soit $f=Z^3-XY$. Notons $h_1(t)$, $h_2(t)$ et $h_3(t)$ les arcs d\'efinis par $h_1(t)=(t^5,\,t^7,\,t^4)$, $h_2(t)=(t^6,\,t^6,\,t^4)$ et $h_3(t)=(t^7,\,t^5,\,t^4)$. Modulo $t^5$, ces trois  arcs ont la m\^eme troncation mais les ordres des diff\'erentes coordonn\'ees ne co\"incident pas. C'est-\`a-dire que $\Chi_{4,\,5,\,7}\cap\Chi_{4,\,6,\,6}\cap \Chi_{4,\,7,\,5}\neq\emptyset$.
\end{exem}
N\'enamoins nous pouvons \'enoncer le r\'esultat suivant :
\begin{lemm}
Nous avons 
$$\Chi_{p,\,l_1,...,\,l_m}=\Chi_{p,\,l'_1,...,\,l'_m}\Longleftrightarrow \Chi_{p,\,l_1,...,\,l_m}\cap\Chi_{p,\,l'_1,...,\,l'_m}\neq\emptyset$$
$$\qquad\qquad\qquad\quad\Longleftrightarrow \left\{\begin{array}{c}l_i\leq p\Rightarrow l_i=l'_i\\
l_i>p \Leftrightarrow l'_i>p\\
b_{r_{\underline{l}}}(\underline{l}-\underline{l}')=0\end{array}\right.$$
o\`u $r_{\underline{l}}$ est l'entier $k_i$ si $l_i>p$ (voir remarque \ref{bonplan}).
\end{lemm}
\begin{proof}
Il est clair que
$$\Chi_{p,\,l_1,...,\,l_m}=\Chi_{p,\,l'_1,...,\,l'_m}\Longrightarrow \Chi_{p,\,l_1,...,\,l_m}\cap\Chi_{p,\,l'_1,...,\,l'_m}\neq\emptyset.$$
L'implication
$$\Chi_{p,\,l_1,...,\,l_m}\cap\Chi_{p,\,l'_1,...,\,l'_m}\neq\emptyset \Longrightarrow \left\{\begin{array}{c}b_{r_{\underline{l}}}(\underline{l}-\underline{l}')=0,\\
l_i\leq p\Rightarrow l_i=l'_i\\
l_i>p \Leftrightarrow l'_i>p\end{array}\right.$$
d\'ecoule de la remarque faite \`a la fin de la preuve de la proposition \ref{arctr}.\\
Fixons $(l_1,...,\,l_m)\in\Ker\,M$ et $(l'_1,...,\,l'_m)\in\Ker\,M$.  Supposons que nous ayons $b_{r_{\underline{l}}}(\underline{l}-\underline{l}')=0$, $l_i\leq p\Rightarrow l_i=l'_i$ et que $l_i>p \Leftrightarrow l'_i>p$.  Supposons qu'il existe un entier $i$ pour lequel $l_i>p$. Dans le cas contraire les \'equivalences sont triviales. Soient $h\in\Chi_{p,\,l_1,...,\,l_m}$ et $\ovl{h}$ un arc \'egal \`a $h$ modulo $(t)^{p+1}$ et dont chaque composante $x_i$ est d'ordre $l_i$ pour $1\leq i\leq m$. Nous avons donc $\ovl{h}(t)=(x_1(t),...,\,x_m(t),\,z(t))$ et nous pouvons supposer, quitte \`a changer les variables, que $x_1(t)=...=x_k(t)=0$ modulo $(t)^{p+1}$, et que $x_{k+1}(t)$,..., $x_m(t)$ sont non nuls modulo $(t)^{p+1}$ (c'est-\`a-dire que $l_1,...,\,l_k>p$ et $l_{k+1},...,\,l_m\leq p$). Nous allons chercher des $x_i'(t)$, avec $\ord(x'_i(t))=l'_i$ pour $1\leq i\leq m$, tels que 
$$\xi(x^{1/n}(t))=\xi(x'^{1/n}(t)).$$
Pour cela, il nous suffit de poser $x'_i(t)=x_i(t)$ pour tout $i$ compris entre $k+1$ et $m$, et de voir que cela revient alors \`a trouver des $x'_i(t)$ tels que 
\begin{equation}\label{eqa}\xi_{k_i}(x^{1/n}(t))+\xi_{k_i+1}(x^{1/n}(t))+\cdots=\xi_{k_i}(x'^{1/n}(t))+\cdots.\end{equation}
Il suffit pour cela de trouver $x'_{i,l'_i}$ pour $k+1\leq i\leq m$ tels que $\prod_ix_{i,l_i}^{a_i(k_i)}=\prod_ix'^{a_i(k_i)}_{i,l'_i}$ ce qui est toujours possible. Ensuite il suffit de choisir les $x'_{i,k}$ pour $k>l'_i$ ce qui est toujours possible car les \'equations en les $x'_{i,k}$, qui d\'ecoulent de l'annulation des diff\'erents termes de l'equation \ref{eqa}, forment un syst\`eme triangulaire. Donc $h\in\Chi_{p,\,l'_1,...,\,l'_m}$. \\
\end{proof}
\begin{exem}
Si pour tous $i$ et $j$ tels que $i\neq j$ nous avons $a_i(1)+a_j(1)\geq 1$, alors l'union (\ref{union}) pr\'ec\'edente est disjointe.
\end{exem}
Nous allons alors noter : 
$$N_{p,\,l_1,...,\,l_m}:=Card\left\{\Ker(\,b_{p,\,l_1,...,\,l_m})\bigcap\prod_{i\in I_{p,\,\underline{l}}}\left(\Z\,\cap]-\infty,\,l_i-p[\right)\right\}$$
o\`u $I_{p,\,\underline{l}}=\{i_1,...,\,i_q\}$ est l'ensemble des indices $i$ pour lesquels $l_i>p$ et $b_{p,\,l_1,...,\,l_m}$ est l'application lin\'eaire
$$b_{p,\,l_1,...,\,l_m}\,:\ \Z^{q}\lgw \Z$$
$$ (u_{1},...,\,u_{q})\lgm na_{i_1}(r_{\underline{l}})u_1+\cdots+na_{i_q}(r_{\underline{l}})u_q$$
D'apr\`es le lemme pr\'ec\'edent, nous voyons que $\Chi_{p,\,l_1,...,\,l_m}=\Chi_{p,\,l'_1,...,\,l'_m}$ si seulement si $(l_1-l'_1,...,\,l_m-l'_m)\in \Ker(\,b_{p,\,l_1,...,\,l_m})\bigcap\prod_{i\in I_{p,\,\underline{l}}}\left(\Z\,\cap]-\infty,\,l_i-p[\right)$.\\
La mesure motivique des $p$-troncations d'arcs qui ne se rel\`event pas en arcs dont l'une des coordonn\'ees $x_i$ est nulle est alors \'egale \`a :
\begin{equation}\left[\bigcup_{(l_1,...,\,l_m)\in D(m)_p}\Chi_{p,\,l_1,...,\,l_m}\right]=\sum_{(l_1,...,\,l_m)\in D(m)_p}\frac{[\Chi_{p,\,l_1,...,\,l_m}]}{N_{p,\,l_1,...,\,l_m}}\end{equation}\\
\begin{exem}
Dans le cas $Z^3=XY$, nous avons 
$$N_{4,\,5,\,7}=Card\left\{\Ker\,(1,\,1)\bigcap\left(\left(\Z\,\cap]-\infty,\,1[\right)\times\left(\Z\,\cap]-\infty,\,3[\right)\right)\right\}.$$
Cet ensemble est form\'e de trois points : le point $(0,\,0)$ correspond \`a $\Chi_{4,\,5,\,7}$, le point $(1,\,-1)$ correspond \`a $\Chi_{4,\,6,\,6}$ et le point $(2,\,-2)$ correspond \`a $\Chi_{4,\,7,\,5}$. Et donc le cardinal de cet ensemble vaut 3.

\end{exem}

\section{\'Etude des troncations d'arcs pour lesquels un des $x_i(t)$ est nul}
Ce terme correspond aux troncations d'arcs pour lesquels un des $x_i(t)$ est nul.\\
Nous allons exprimer le terme $[X(a(1),...,\,a(g))_{p,cn}]$ en fonction de diff\'erents termes de la forme $[\pi_p(X(b(1),...,\,b(i))_{\infty})]$ pour $i<g$.\\
Rappelons que pour tout $i$ compris entre 1 et $m$, $k_i$ est le plus petit entier qui v\'erifie $a_{i}(k_i)\neq 0$ et que nous avons
$$1=k_1=k_2=\cdots=k_{i_0}<k_{i_0+1}\leq k_{i_0+2}\leq\cdots\leq k_m\leq g$$
Soit $i\in\{1,...,\,m\}$. Que vaut $[\Chi_{p,\,l_1,...,\,l_m}]$ o\`u $l_i=+\infty$ et $l_j<+\infty$ pour $j<i$ ?\\
 Si $i>i_0$ les $x_j$ pour $j>i$ peuvent \^etre choisis quelconques  et 
$$[\Chi_{p,\,l_1,...,\,l_m}]=\DL^{p(m-i)}[\pi_p(X(a^i(1),...,\,a^i(k_i-1))_{\infty})]$$
avec $a^i(j)=(a_1(j),...,\,a_{k_i-1}(j))$ pour $1\leq j\leq k_i-1$. Le germe de vari\'et\'e $X(a^i(1),...,\,a^i(k_i-1))$ n'est peut-\^etre pas r\'eduit mais il reste irr\'eductible car ses exposants caract\'eristiques sont conjugu\'es.\\
Si $i\leq i_0$, les $x_j$ pour $j\neq i$  peuvent \^etre choisis quelconques mais les  $x_j$ pour $j<i$ doivent \^etre non nuls donc $[\Chi_{p,\,l_1,...,\,l_m}]=\DL^{p(m-i)}\left(\DL^p-1\right)^{i-1}$.\\

Nous avons donc
\begin{equation}\label{2eme terme}\begin{split}[X(a(1),...,\,a(g))_{p,cn}]=&\sum_{i=1}^{i_0}\DL^{p(m-i)}\left(\DL^{p}-1\right)^{i-1}\\
+\sum_{i=i_0+1}^m\DL^{p(m-i)}&[\pi_p(X(a^i(1),...,\,a^i(k_i-1))_{\infty})],\\
\end{split}\end{equation}
et 
\begin{multline}P_{geom}^{tore^c}(a(1),...,\,a(m))(T)=\sum_{i=0}^{i_0-1}\sum_{j=0}^i\frac{C_{i}^{i-j}(-1)^{j}}{1-\DL^{m-j-1}T}\\
+\sum_{i=i_0+1}^mP_{geom,X(a^i(1),...,\,a^i(k_i-1)),0}(\DL^{m-i}T).
\end{multline}
En suivant le m\^eme raisonnement nous obtenons
\begin{equation}\label{2eme terme'}\begin{split}\Chi_c(\varphi_{p,cn}(a(1),...,\,a(m)))=&\sum_{i=1}^{i_0}\DL^{p(m-i)}\left(\DL^{p}-1\right)^{i-1}\\
+\sum_{i=i_0+1}^m\DL^{p(m-i)}&\Chi_c(\varphi_p(a^i(1),...,\,a^i(k_i-1))),\\
\end{split}\end{equation}
et 
\begin{multline}P_{arit}^{tore^c}(a(1),...,\,a(m))(T)=\sum_{i=0}^{i_0-1}\sum_{j=0}^i\frac{C_{i}^{i-j}(-1)^{j}}{1-\DL^{m-j-1}T}\\
+\sum_{i=i_0+1}^mP_{arit,X(a^i(1),...,\,a^i(k_i-1)),0}(\DL^{m-i}T).
\end{multline}

\section{R\'esultat principal}
Nous pouvons alors donner la forme g\'en\'erale des deux s\'eries \'etudi\'ees :
\begin{theo}\label{series}
Soit $(X,\,0)$  un germe d'hypersurface irr\'eductible quasi-ordinaire d'exposants caract\'eristiques $a(1),...,\,a(g)$. Nous avons alors les relations de r\'ecurrence suivantes
\begin{equation}\begin{split}
P_{geom,X(a(1),...,\,a(g)),0}(T)=(\DL-1)^m\sum_{i=1}^m\sum_{p\geq0}\DL^{pm}T^p\\
\times\sum_{(l_1,...,\,l_m)\in D_{0,i}(m)_p}\DL^{-\sum_{j=1}^ml_j+\max\{l_i-b_{k_i}(\underline{l}),\,0\}}+\\
+\sum_{q=1}^{m}\sum_{i=1}^m(\DL-1)^{m-q+1}\sum_{p\geq 0}(\DL^{m-q+1}T)^p\times\sum_{(l_1,...,\,l_m)\in D_{q,i}(m)_p}\frac{\DL^{-\sum_{j\notin I_{p,\,\underline{l}}}l_j-b_{k_i}(\underline{l})}}{N_{p,\,l_1,...,\,l_m}}+\\
+\sum_{i=0}^{i_0-1}\sum_{j=0}^i\frac{C_{i}^{i-j}(-1)^{j}}{1-\DL^{m-j-1}T}
+\sum_{i=i_0+1}^mP_{geom,X(a^i(1),...,\,a^i(k_i-1)),0}(\DL^{m-i}T)
\end{split}
\end{equation}
\begin{equation}\begin{split}
P_{arit,X(a(1),...,\,a(g)),0}(T)=(\DL-1)^m\sum_{k=0}^g\frac{n_{k_l}}{n}\sum_{i=1}^m\sum_{p\geq0}\DL^{pm}T^p\\
\times\sum_{(l_1,...,\,l_m)\in D_{0,i}(m)_p}\DL^{-\sum_{j=1}^ml_j+\max\{l_i-b_{k_i}(\underline{l}),\,0\}}+\\
+\sum_{k=0}^g\sum_{q=1}^{m}\sum_{i=1}^m\sum_{\begin{array}{c}I\subset\{1,...,\,m\}\\
i\in I, \#I=q\end{array}}(\DL-1)^{m-q+1}\sum_{p\geq 0}(\DL^{m-q+1}T)^p\times\\
\times\sum_{(l_1,...,\,l_m)\in D_{q,i, I}(m)_p}\frac{1}{n(i,I_{p,\,\underline{l}})}\frac{\DL^{-\sum_{j\notin I_{p,\,\underline{l}}}l_j-b_{k_i}(\underline{l})}}{N_{p,\,l_1,...,\,l_m}}+\\
+\sum_{i=0}^{i_0-1}\sum_{j=0}^i\frac{C_{i}^{i-j}(-1)^{j}}{1-\DL^{m-j-1}T}
+\sum_{i=i_0+1}^mP_{geom,X(a^i(1),...,\,a^i(k_i-1)),0}(\DL^{m-i}T)
\end{split}\end{equation}
o\`u 
$$D(m)_p:=\left\{(l_1,...,\,l_m)\in\Ker\,M\,/\, 0<l_i\leq p\right\}$$
$$\bigcup\left(\bigcup_{i=1}^m\left\{(l_1,...,\,l_m)\in\Ker\,M\,/\,l_i>p,\,l_j\leq p,\,j=i+1,...,m,\,b_{k_i}(\underline{l})\leq p\right\}\cap(\N^*)^m\right)\, ,$$
$I_{p,\,\underline{l}}$ est l'ensemble des indices $j$ pour lesquels $l_j>p$,
$$D_{0,i}(m)_p:=\left\{\underline{l}\in D(m)_p\ \backslash\ l_j\leq p \ \ \forall j,\ l_i-b_{k_i}(\underline{l})\geq l_j-b_{k_j}(\underline{l}) \ \ \forall j\neq i \right.$$
$$\left.\text{ et }l_i-b_{k_i}(\underline{l})> l_j-b_{k_j}(\underline{l})\  \ \forall j<i\right\},$$
$$D_{q,i}(m)_p:=\left\{\underline{l}\in D(m)_p\ \backslash\ l_i> p, \ 
Card(I_{p,\underline{l}})=q,\ l_i-b_{k_i}(\underline{l})\geq l_j-b_{k_j}(\underline{l}) \ \ \forall j\neq i \right.$$
$$\left.\text{ et }l_i-b_{k_i}(\underline{l})> l_j-b_{k_j}(\underline{l})\  \ \forall j<i\right\}.$$
$$D_{q,i,I}(m)_p:=\left\{\underline{l}\in D(m)_{q,i}(m)_p\,\backslash\, I_{p,\,\underline{l}}=I\right\},$$
$$\text{et }\ \begin{array}{ccc} M\ :\ \Z^m\ & \lgw & \left(\frac{\Z}{n\Z}\right)^g\\
(l_1,...,\,l_m)& \lgm & (n\sum a_i(1)l_i,...,\, n\sum a_i(g)l_i)\end{array}$$\\
\end{theo}

\section{Cas o\`u les $a_i(1)\geq1$ pour tout $i$}
Dans ce cas beaucoup de choses peuvent \^etre simplifi\'ees, et nous pouvons donner une expression explicite de ces deux s\'eries. Tous les $k_j$ sont \'egaux \`a 1, et $b_1(\underline{l})\geq l_j$ pour tout $j$, si $l_j>0$ pour tout $j$. En particulier, si $l_i>p$, alors $b_{k_i}(\underline{l})>p$ et seul le cas C1 est \`a consid\'erer ici. D'autre part $e=0$. Nous avons
$$D(m)_p=\left\{(l_1,...,\,l_m)\in\Ker\,M\,/\, 0<l_i\leq p\right\}.$$
Notons 
$$\E_p:=\sum_{\begin{array}{c}(l_1,...,\,l_m)\in \Ker\,M\\
0<l_i\leq p\end{array}}\DL^{-\sum_{i=1}^pl_i}$$
$$\text{et }\F_{k,p}:=\sum_{\begin{array}{c}(l_1,...,\,l_m)\in \Ker\,M\\
0<l_i\leq p\\
b_k(\underline{l})\leq p\end{array}}\DL^{-\sum_{i=1}^pl_i}\ .$$
Nous avons alors
$$P_{geom}^{tore}(a(1),...,\,a(m))(T)=(\DL-1)^m\sum_{p\geq 0}\E_p(\DL^{m}T)^p$$
et
\begin{equation}\label{A}
\begin{split}
P_{arit}^{tore}(a(1),...,\,a(m))(T)& =\sum_{k=1}^{g-1}(\DL-1)^m\sum_{p\geq 0}(\DL^{m}T)^p\left(\F_{k,\,p}-\F_{k+1,\,p}\right)\\
+(\DL-1)^m&\sum_{p\geq 0}(\DL^{m}T)^p\F_{g,\,p}+(\DL-1)^m\sum_{p\geq 0}(\E_p-\F_{1,\,p})(\DL^{m}T)^p\\
=(\DL-1)^m\sum_{p\geq 0}&\E_p(\DL^{m}T)^p
+(\DL-1)^m\sum_{k=1}^{g}\left(\frac{n_k-n_{k-1}}{n}\right)\sum_{p\geq 0}\F_{k,p}(\DL^mT)^p \ .
\end{split}\end{equation}
Pour calculer ces s\'eries, notons 
$$\E_{\underline{p}}=\E_{p_1,...,\,p_m}=\sum_{\begin{array}{c}(l_1,...,\,l_m)\in \Ker\,M\\
0<l_i\leq p_i\end{array}}\DL^{-\sum_{i=1}^pl_i}.$$
Nous avons alors le lemme suivant

\begin{lemm}\label{E}
Nous avons pour tout $n\geq p_k> 0$, $p>0$, $k\geq 0$ :
\begin{enumerate}
\item$$\E_p=\E_{p,...,\,p}\ ,$$
\item$$\E_{p_1,...,\,p_j+kn,...,\,p_m}=\E_{p_1,...,\,kn,...,\,p_m}+\DL^{-kn}\E_{p_1,...,\,p_j,...,\,p_m}\ ,$$
\item$$\E_{p_1,...,\,kn,...,\,p_m}=\frac{1-\DL^{-kn}}{1-\DL^{-n}}\E_{p_1,...,\,n,...,\,p_m}\ ,$$
\item
$$\E_{p_1,...,\,p_j+kn,...,\,p_m}=\frac{1-\DL^{-kn}}{1-\DL^{-n}}\E_{p_1,...,\,n,...,\,p_m}+\DL^{-kn}\E_{p_1,...,\,p_j,...,\,p_m}\ .$$
\end{enumerate}
Notons $\mathcal{S}_m:=\left\{s : \{1,...,\,m\}\lgw \{0,1\}^m\right\}$. Pour tout $s\in\mathcal{S}_m$, nous noterons $|s|$ le nombre d'\'el\'ements de $\{1,...,\,m\}$ d'image $1$ par $s$. Nous avons alors
$$\E_{\underline{p}+n\underline{k}}=\sum_{s\in\mathcal{S}_m}\DL^{\sum_i-k_in(1-s(i))}\E_{p_1(1-s(1))+k_1ns(1),...,\,p_m(1-s(m))+k_mns(m)}\ .$$
En notant $\E_{p,s}=\E_{p(1-s(1))+ns(1),...,\,p(1-s(m))+ns(m)}$, nous avons
$$\E_{p+nk}=\sum_{s\in\mathcal{S}_m}\left(\frac{\DL^{-kn}-1}{\DL^{-n}-1}\right)^{|s|}\DL^{-kn(m-|s|)}\E_{p,s}\ .$$
\end{lemm}

\begin{proof} Le 1. d\'ecoule des d\'efinitions.\\
Il faut tout d'abord remarquer que nous avons $M(l_1,...,\,l_m)=0$ si et seulement si $M(l_1,...,\,l_i+ln,...,\,l_m)=0$ pour tout $i$ et tout $l$ entier. Le 2. d\'ecoule  alors de l'\'egalit\'e suivante :$$\E_{p_1,...,\,p_j+kn,...,\,p_m}=\sum_{\begin{array}{c}(l_1,...,\,l_m)\in \Ker\,M\\
0<l_i\leq p_i\ i\neq j\\
0<l_j\leq p_j+kn\end{array}}\DL^{\sum_{i=1}^p-l_i}$$
$$=\sum_{\begin{array}{c}(l_1,...,\,l_m)\in \Ker\,M\\
0<l_i\leq p_i\ i\neq j\\
0<l_j\leq kn\end{array}}\DL^{\sum_{i=1}^p-l_i}+\DL^{-kn}\sum_{\begin{array}{c}(l_1,...,\,l_m)\in \Ker\,M\\
0<l_i\leq p_i\ \end{array}}\DL^{\sum_{i=1}^p-l_i}\ .$$
Le 3. s'en d\'eduit par r\'ecurrence et le 4. en compilant 2. et 3.\\
\end{proof}

\begin{coro}\label{E'}
Nous avons 
$$\sum_{p\geq 0}\E_p(\DL^{m}T)^p=\sum_{p=1}^n\sum_{s\in\mathcal{S}_m}(\DL^mT)^p\frac{\E_{p,s}}{(\DL^{-n}-1)^{|s|}}\sum_{j=0}^{|s|}\frac{C^j_{|s|}(-1)^j}{1-(\DL^jT)^n}\ . $$

\end{coro}
\begin{proof} Nous avons
$$\sum_{p\geq 0}\E_p(\DL^{m}T)^p=\sum_{p=1}^n\sum_{k\geq0} (\DL^mT)^{kn+p}\E_{p+kn}$$
$$=\sum_{p=1}^n(\DL^mT)^p\sum_{s\in\mathcal{S}_m}\sum_{k\geq0}(\DL^mT)^{nk}\left(\frac{\DL^{-kn}-1}{\DL^{-n}-1}\right)^{|s|}\DL^{-kn(m-|s|)}\E_{p,s}\ .$$
Le r\'esultat est alors direct en d\'eveloppant $\left(\frac{\DL^{-kn}-1}{\DL^{-n}-1}\right)^{|s|}$ et sommant sur $k$.\\\end{proof}

Pour calculer les s\'eries o\`u apparaissent les $\F_{k,p}$ nous pouvons remarquer comme tous les $a_i(k)\geq 1$ que 
$$\F_{k,p}=\sum_{\begin{array}{c}(l_1,...,\,l_m)\in \Ker\,M\\
0<l_i\\
b_k(\underline{l})\leq p\end{array}}\DL^{-\sum_{i=1}^pl_i}\ .$$
Donc nous avons
$$\sum_{\begin{array}{c}p\geq 0\end{array}}\F_{k,p}(\DL^mT)^p=\sum_{\begin{array}{c}p\geq 0\end{array}}\sum_{\begin{array}{c}(l_1,...,\,l_m)\in \Ker\,M\\
0<l_i\\
b_k(\underline{l})\leq p\end{array}}\DL^{-\sum_{i=1}^pl_i}(\DL^mT)^p$$
$$=\sum_{\begin{array}{c}(l_1,...,\,l_m)\in \Ker\,M\\
0<l_i\\
\end{array}}\sum_{\begin{array}{c}p\geq b_k(\underline{l})\end{array}}\DL^{-\sum_{i=1}^pl_i}(\DL^mT)^p$$
$$=\sum_{\begin{array}{c}(l_1,...,\,l_m)\in \Ker\,M\\
0<l_i\\
\end{array}}\frac{\DL^{-\sum_{i=1}^pl_i}(\DL^mT)^{b_k(\underline{l})}}{1-\DL^mT}$$
$$=\sum_{\begin{array}{c} r_i\geq 0\end{array}}\sum_{\begin{array}{c}(l_1,...,\,l_m)\in \Ker\,M\\
0<l_i\leq n\\
\end{array}}\frac{\DL^{-\sum_{i=1}^pl_i}(\DL^mT)^{b_k(\underline{l})}}{1-\DL^mT}\DL^{-n\sum_{i=1}^p r_i}(\DL^mT)^{nb_k(\underline{r})}\ .$$
Car $(l_1+nr_1,...,\,l_m+nr_m)\in \Ker\,M\Longleftrightarrow (l_1,...,\,l_m)\in \Ker\,M$. Nous avons donc
$$\sum_{\begin{array}{c}p\geq 0\end{array}}\F_{k,p}(\DL^mT)^p=$$
$$=\sum_{\begin{array}{c}(l_1,...,\,l_m)\in \Ker\,M\\
0<l_i\leq n\\
\end{array}}\frac{\DL^{-\sum_{i=1}^pl_i}(\DL^mT)^{b_k(\underline{l})}}{1-\DL^mT}\prod_{i=1}^m \frac{1}{1-\DL^{n(ma_i(k)-1)}T^{na_i(k)}}\ .$$
Nous pouvons alors en d\'eduire le
\begin{theo}\label{exp}
Soit $(X,0)$ est un germe irr\'eductible \`a singularit\'e quasi-ordinaire tel  que $a_i(1)\geq 1$ pour tout $i\in\{1,...,m\}$. Alors nous avons
\begin{equation}\begin{split}P_{geom,X,0}(T)&=\sum_{i=0}^{m-1}\sum_{j=0}^i\frac{C_{i}^{i-j}(-1)^{j}}{1-\DL^{m-j-1}T}\\
+(\DL-1)^m&\sum_{p=1}^n\sum_{s\in\mathcal{S}_m}(\DL^mT)^p\frac{\E_{p,s}}{(\DL^{-n}-1)^{|s|}}\sum_{j=0}^{|s|}\frac{C^j_{|s|}(-1)^j}{1-(\DL^jT)^n}\ ,
\end{split}\end{equation}

\begin{equation}\begin{split}P_{arit,X,0}(T)&=\sum_{i=0}^{m-1}\sum_{j=0}^i\frac{C_{i}^{i-j}(-1)^{j}}{1-\DL^{m-j-1}T}\\
+\frac{1}{n}(\DL-1)^m&\sum_{p=1}^n\sum_{s\in\mathcal{S}_m}(\DL^mT)^p\frac{\E_{p,s}}{(\DL^{-n}-1)^{|s|}}\sum_{j=0}^{|s|}\frac{C^j_{|s|}(-1)^j}{1-(\DL^jT)^n}\\
+(\DL-1)^m&\sum_{k=1}^{g}\left(\frac{n_k-n_{k-1}}{n}\right)\frac{\H_{k}}{1-\DL^mT}\prod_{i=1}^m \frac{1}{1-\DL^{n(ma_i(k)-1)}T^{na_i(k)}}\\
\end{split}\end{equation}
o\`u $\E_{p,s}$ et $\E_{\underline{p},s}$ sont d\'efinis lemme \ref{E} et 
$$\H_{k}:=\sum_{\begin{array}{c}(l_1,...,\,l_m)\in \Ker\,M\\
0<l_i\leq n\\
\end{array}}\DL^{-\sum_{i=1}^pl_i}(\DL^mT)^{b_k(\underline{l})}$$

\end{theo}
\subsection{Exemples}
\subsubsection{Courbes planes}
Consid\'erons un germe de courbe plane de caract\'eristique $(\b_0=n,\,\b_1,...,\b_g)$ \cite{Z}. Nous avons alors $n=n$, $m=1$ et $g=g$. De plus
$$l\in\Ker\,M \Longleftrightarrow l\in n\Z\ .$$
Nous obtenons alors
$$P_{geom,X,0}(T)=\frac{1}{1-T}+\frac{\DL-1}{1-\DL T}\frac{T^n}{1-T^n}$$
$$\text{et }P_{arit,X,0}(T)=\frac{1}{1-T}+\frac{\DL-1}{1-\DL T}\sum_{k=0}^{g}\left( \frac{n_k-n_{k-1}}{n}\right)\frac{\DL^{\b_k-n}T^n}{1-\DL^{\b_k-n}T^n}.$$
On retrouve l\`a le r\'esultat de Denef et Loeser \cite{D-L2}.
\subsubsection{L'hypersurface d\'efinie par $Z^2-X^3Y^3=0$}
Soit l'hypersurface de $\C^3$ d\'efinie par $Z^2-X^3Y^3=0$. La s\'erie fractionnaire associ\'ee est $\xi=X^{3/2}Y^{3/2}$. Dans ce cas nous avons $n=2$, $m=2$, $g=1$, $n_1=2$  et

$$\begin{array}{ccc} M\ :\ \Z^2\ & \lgw & \frac{\Z}{2\Z}\\
(l_1,\,l_2)& \lgm & 3l_1+3l_2\end{array}$$
Nous avons 
$$(l_1,\,l_2)\in\Ker\ M\Longleftrightarrow l_1+l_2\in  \frac{\Z}{2\Z}\ .$$
Soit $p\in\{1,\,2\}$, nous avons 
$$\E_{1,(0,0)}=\E_{1,(1,0)}=\E_{1,(0,1)}=\DL^{-2},$$
$$\E_{2,(0,0)}=\E_{2,(0,1)}=\E_{2,(1,0)}=\E_{2,(1,1)}=\E_{1,(1,1)}=\DL^{-2}+\DL^{-4}\ .$$
D'o\`u

$$(1+\DL)^2P_{geom,X,0}(T)=\frac{-2\DL}{1-T}+\frac{(1+\DL)^2}{1-\DL T}-\frac{(1-\DL)^2}{1+\DL T}+\frac{1+\DL^2}{1-\DL^2T}\ .$$
Dans ce cas, la s\'erie g\'eom\'etrique a 4 p\^oles en $T$ : $1$,  $\DL^{-1}$, $-\DL^{-1}$ et $\DL^{-2}$.

\end{document}